  \documentclass[12pt]{amsart}
  \usepackage{graphicx}
  \usepackage{latexsym} 
  \usepackage[all]{xy}  
  \usepackage{amsfonts} 
  \usepackage{amsthm}   
  \usepackage{amsmath}  
  \usepackage{amssymb}  
  \usepackage{enumerate}
  \usepackage[2emode]{psfrag}
  \usepackage{hyperref}
  \usepackage{color}  
  \newcommand{\cC}{{C}}  

\numberwithin{equation}{section}

\DeclareMathOperator{\newdiamond}{\raisebox{0.2ex}{\scalebox{0.75}{\ensuremath{\lozenge}}}}
\DeclareMathOperator{\newblackdiamond}{\raisebox{0.2ex}{\scalebox{0.75}{\ensuremath{\blacklozenge}}}}
\DeclareMathOperator{\newdiamondd}{\raisebox{0.2ex}{\scalebox{0.5}{\ensuremath{\lozenge}}}}
\DeclareMathOperator{\newblackdiamondd}{\raisebox{0.2ex}{\scalebox{0.5}{\ensuremath{\blacklozenge}}}}

\DeclareMathOperator{\newtriangle}{\raisebox{0ex}{\scalebox{0.75}{\ensuremath{\triangledown}}}}
\DeclareMathOperator{\newblacktriangle}{\raisebox{0ex}{\scalebox{0.75}{\ensuremath{\blacktriangledown}}}}
\DeclareMathOperator{\newtriangled}{\raisebox{0.1ex}{\scalebox{0.5}{\ensuremath{\triangledown}}}}
\DeclareMathOperator{\newblacktriangled}{\raisebox{0.1ex}{\scalebox{0.5}{\ensuremath{\blacktriangledown}}}}

\DeclareMathOperator{\newtriangleup}{\raisebox{0ex}{\scalebox{0.75}{\ensuremath{\vartriangle}}}}
\DeclareMathOperator{\newblacktriangleup}{\raisebox{0ex}{\scalebox{0.75}{\ensuremath{\blacktriangle}}}}

\def\cC{\mathcal C}
\def\cD{\mathcal D}
\def\cM{{\mathcal M}} 
\def\can{{\rm can}} 
\def\beq{\begin{equation}} 
\def\eeq{\end{equation}}

\def\ot{{\otimes}}

\newcounter{zlist} 
  \newenvironment{zlist}{\begin{list}{(\arabic{zlist})}{ 
  \usecounter{zlist}\leftmargin2.5em\labelwidth2em\labelsep0.5em 
  \topsep0.6ex
  \parsep0.3ex plus0.2ex minus0.1ex}}{\end{list}}

\newcounter{blist} 
  \newenvironment{blist}{\begin{list}{(\alph{blist})}{ 
  \usecounter{blist}\leftmargin2.5em\labelwidth2em\labelsep0.5em 
  \topsep0.6ex 
  \parsep0.3ex plus0.2ex minus0.1ex}}{\end{list}}

\def\stac#1{\raise-.2cm\hbox{$\stackrel{\displaystyle\,\otimes\,}{\scriptscriptstyle{#1}}$}}
\def\sstac#1{\otimes_{#1}}
\def\ddiamonda{{\raise1.2pt\hbox{$\lozenge$}}\!\!\!\!{\raise-1.2pt\hbox{$\lozenge$}}}
\def\blackddiamonda{{\raise1.2pt\hbox{$\blacklozenge$}}\!\!\!\!{\raise-1.2pt\hbox{$\blacklozenge$}}}

\def\ddiamond{{{\raisebox{0.3ex}{\resizebox{1.7ex}{!}{\ensuremath{\,\ddiamonda\,}}}}}}
\def\blackddiamond{{{\raisebox{0.3ex}{\resizebox{1.7ex}{!}{\ensuremath{\,\blackddiamonda\,}}}}}}

\newtheorem{proposition}{Proposition}[section] 
  \newtheorem{lemma}[proposition]{Lemma} 
   
  \newtheorem{corollary}[proposition]{Corollary} 
  \newtheorem{theorem}[proposition]{Theorem} 

\theoremstyle{definition} 
  \newtheorem{definition}[proposition]{Definition}

  \theoremstyle{remark} 
  \newtheorem{remark}[proposition]{Remark}

  \newcounter{c} 
  \renewcommand{\[}{\setcounter{c}{1}$$} 
  \newcommand{\etyk}[1]{\vspace{-7.4mm}$$\begin{equation}\Label{#1} 
  \addtocounter{c}{1}} 
  \renewcommand{\]}{\ifnum \value{c}=1 $$\else \end{equation}\fi} 
\def\squar{\raisebox{0.05ex}{\tiny{\ensuremath{\square}}}}
\def\blacksquar{\raisebox{0.05ex}{\tiny{\ensuremath{\blacksquare}}}}

\def\tildej{\widetilde{\jmath}}
\def\tildeJ{\widetilde{J}}

\newcommand{\thlabel}[1]{\label{th:#1}}
\newcommand{\thref}[1]{Theorem~\ref{th:#1}}
\newcommand{\selabel}[1]{\label{se:#1}}
\newcommand{\seref}[1]{Section~\ref{se:#1}}
\newcommand{\lelabel}[1]{\label{le:#1}}
\newcommand{\leref}[1]{Lemma~\ref{le:#1}}
\newcommand{\prlabel}[1]{\label{pr:#1}}
\newcommand{\prref}[1]{Proposition~\ref{pr:#1}}
\newcommand{\colabel}[1]{\label{co:#1}}
\newcommand{\coref}[1]{Corollary~\ref{co:#1}}
\newcommand{\relabel}[1]{\label{re:#1}}
\newcommand{\reref}[1]{Remark~\ref{re:#1}}

\newcommand{\eqlabel}[1]{\label{eq:#1}}
\newcommand{\equref}[1]{(\ref{eq:#1})}

\newcommand{\Hom}{{\rm Hom}}
\newcommand{\End}{{\rm End}}

\def\ul{\underline}

\def\ot{\otimes}

\newcommand{\Cc}{\mathcal{C}}
\newcommand{\Dd}{\mathcal{D}}

\newcommand{\Mm}{\mathcal{M}}

\def\*C{{}^*\hspace*{-1pt}{\Cc}}

\def\text#1{{\rm {\rm #1}}}

\begin{document} 

 \title{Morita theory for Coring extensions and Cleft bicomodules} 
 \author{Gabriella B\"ohm}
 \address{Research Institute for Particle and Nuclear Physics, Budapest, 
 \newline\indent H-1525
 Budapest 114, P.O.B.\ 49, Hungary}
  \email{G.Bohm@rmki.kfki.hu}
  \author{Joost Vercruysse}    
  \address{Vrije Universiteit Brussel VUB, Pleinlaan 2, B-1050,
  Brussel, Belgium} 
 \email{joost.vercruysse@vub.ac.be}   
    \date{April 2006} 
  \subjclass{16W30} 
  \begin{abstract}
 A {\em Morita context} is constructed for any comodule of a coring and,
 more generally, for an $L$-$\cC$ bicomodule $\Sigma$ for a 
 pure coring extension
 $(\cD:L)$ of $(\cC:A)$. It is related to a 2-object subcategory
 of the category of $k$-linear functors $\Mm^\Cc\to\Mm^\Dd$. Strictness of
 the Morita 
 context is shown to imply the Galois property of $\Sigma$
 as a $\cC$-comodule and a Weak Structure Theorem.
 Sufficient conditions are found also for a Strong Structure
 Theorem to hold.

 {\em Cleft} property of an $L$-$\cC$ bicomodule $\Sigma$ --
 implying strictness of the associated Morita context -- is
 introduced. It is shown to be equivalent to being a {\em Galois}
 $\cC$-comodule and isomorphic 
 to $\End^\cC(\Sigma)\otimes_{L} \cD$, in the category of left
 modules for the ring $\End^\cC(\Sigma)$ and right comodules
 for the coring $\cD$, i.e. satisfying the {\em normal basis} property.

 Algebra extensions, that are cleft extensions by a Hopf algebra, a
 coalgebra or a pure
Hopf algebroid, as well as cleft entwining structures
 (over commutative or
 non-commutative base rings) and cleft weak entwining structures, are
 shown to provide examples of cleft bicomodules. 

Cleft extensions by arbitrary Hopf algebroids are described in terms of Morita
contexts that do not necessarily correspond to coring extensions.
\end{abstract}
  \maketitle

\section{Introduction}
The application of {\em Morita theory} in the study of algebra
extensions by a Hopf algebra (i.e. of comodule algebras for a Hopf
algebra) originates in the work \cite{CohFishMont:Morita} of
Cohen, Fischman and Montgomery. In that paper a Morita context has
been associated to a comodule algebra of a Hopf algebra $H$,
under the assumption that $H$ is a finite dimensional algebra over
a field (or a Frobenius algebra over a commutative ring). The
construction has been extended by Doi in \cite{Doi:smmor} to
arbitrary Hopf algebras $H$. Using the observation \cite{Brz:str},
that a comodule algebra of a Hopf algebra can be considered as a
special instance of a base algebra of a coring with a grouplike
element, in \cite{CaVerWang:cleftentw} a Morita context has been
constructed for any coring with a grouplike element. Recall that
the existence of a grouplike element in an $A$-coring $\cC$ is
equivalent to the existence of a (left or right) $\cC$ comodule
structure in $A$. In the paper \cite{CaDeGrVe:comcor} the
particular $\cC$-comodule $A$ in \cite{CaVerWang:cleftentw} was
replaced by any $\cC$-comodule $\Sigma$, which is finitely
generated and projective as an $A$-module. In all of the listed
(more and more general) situations, strictness of the Morita
context was related to properties of the extension (or the coring)
behind.

The first aim of the present paper is to generalize the
construction of a Morita context further. As a first step, in
Section \ref{sec:Morita} we remove the finitely generated
projectivity condition in \cite{CaDeGrVe:comcor} and construct a
Morita context for an arbitrary comodule $\Sigma$ of an
$A$-coring $\cC$. It connects the algebra of $\cC$-comodule
endomorphisms of $\Sigma$
and the $A$-dual algebra of $\cC$.

As a next step, in \seref{contextext} we study pure coring
extensions $(\cD:L)$ of $(\cC:A)$.
By a definition in \cite{Brz:corext}, a coring $\cD$ over a $k$-algebra $L$
is a right extension of another coring $\cC$ over a $k$-algebra $A$ if $\cC$ 
is a $\cC$-$\cD$ bicomodule with the left regular $\cC$-comodule
structure (i.e. via the coproduct $\Delta_\cC$). 
A coring extension $\cD$ of $\cC$ is said to be pure if 
$$
\xymatrix{
M \ar[rr]^-{\varrho}&&
M \ot_A {\mathcal C} \ar@<2pt>[rr]^-{\varrho\ot_A {\mathcal C}}
\ar@<-2pt>[rr]_-{M\ot_A \Delta_{\mathcal C}}&&
M \ot_A {\mathcal C} \ot_A {\mathcal C}
}
$$
is a $\cD\ot_L \cD$-pure equalizer in $\cM_L$,  
for any right $\cC$-comodule $(M,\varrho)$.
By \cite[22.3]{BrzWis:cor} and its Erratum
(see also the arXiv version of \cite[Theorem 2.6]{Brz:corext}), a pure
coring extension induces
a $k$-linear functor $U:\cM^\cC\to \cM^\cD$, such that the
forgetful functor $\cM^\cC\to \cM_k$ factorizes through $U$ and the
forgetful functor $\cM^\cD\to \cM_k$. In this situation, to any
$L$-$\cC$ bicomodule $\Sigma$ (i.e. left $L$-module and right
$\cC$-comodule, with left $L$-linear $\cC$-coaction) we associate
a Morita context. It reduces to the Morita context in Section
\ref{sec:Morita} if the coring $\cD$ is trivial (i.e. equal to the
ground ring $k$).

Note that -- similarly to the identification of $k$-algebras with
$k$-linear categories of a single object -- one can identify
Morita contexts with $k$-linear categories of two objects (cf.
\cite[Remark 3.2 (1)] {BohmBrz:cleft}). Namely, to a $k$-linear
category with two objects $a$ and $b$, associate a Morita context
of the two endomorphism algebras, ${\rm End}(a)$ and ${\rm
End}(b)$. The $k$-modules ${\rm Hom}(a,b)$ and ${\rm Hom}(b,a)$, of 
morphisms between the two objects, are bimodules via the composition 
of morphisms. 
Since the restrictions of the composition of morphisms in the category
to ${\rm Hom}(a,b)\otimes_k {\rm Hom}(b,a)$ and to ${\rm
Hom}(b,a)\otimes_k {\rm Hom}(a,b)$ are balanced in the appropriate
sense, both connecting homomorphisms can be constructed as their
projections. Via this identification, the Morita context in
\seref{contextext} corresponds to the full subcategory of the
category of $k$-linear functors with two objects: the functor
$U:\cM^\cC\to \cM^\cD$, constructed by Brzezi\'nski in
\cite{Brz:corext} for a pure coring extension $(\cD:L)$ of
$(\cC:A)$ on 
one hand, and the functor ${\rm Hom}^\cC(\Sigma,-)\otimes_L
\cD:\cM^\cC\to \cM^\cD$, given in terms of the $L$-$\cC$
bicomodule $\Sigma$, on the other hand.

\seref{cleft} is devoted to a study of {\em cleft} $L$-$\Cc$ bicomodules. As
in \cite{Abuh:cleftentw} and \cite{CaVerWang:cleftentw}, our main tool is
provided by Morita theory. However, as mentioned in the introduction of
\cite{CaVerWang:cleftentw}, using Morita contexts associated to 
comodules for one coring only (say, as in Section \ref{sec:Morita}), it does
not seem to be possible to go beyond a study of cleft entwining 
structures in \cite{Abuh:cleftentw} and \cite{CaVerWang:cleftentw}. On the
other hand, our language of Morita contexts for pure coring extensions, 
developed in \seref{contextext}, gives a very natural framework to introduce
the {cleft property} of an $L$-$\cC$ bicomodule $\Sigma$ for
a pure coring extension $(\cD:L)$ of $(\cC:A)$. As a guiding example, 
consider a Hopf algebra $H$ over a commutative ring $k$ and its
right comodule algebra $A$. Denote the coinvariants of $A$ by $B$.
The $k$-module
$A\otimes_k H$ inherits a left $A$-module structure of $A$ and it
can be equipped with the structure of a right $A$-module in terms
of the right $H$-coaction in $A$.
The $A$-$A$ bimodule $A\otimes_k H$ is an $A$-coring with
coproduct and counit inherited from $H$. The $k$-coring (coalgebra)
$H$ is a pure right extension of the $A$-coring $A\otimes_k H$
and $A$ is 
a $k\,$-$\,(A\sstac k H)$ bicomodule. The Morita context, corresponding to it, 
consists of two subalgebras of the convolution algebra ${\rm
Hom}_k(H,A)$. The two bimodules are given by the sets of right
$H$-comodule maps $H\to A$, w.r.t. the regular $H$-coaction in
$H$, $h\mapsto h_{(1)}\otimes_k h_{(2)}$, and the twisted
coaction, $h\mapsto h_{(2)}\otimes_k S(h_{(1)})$, respectively,
where $S$ denotes the antipode of $H$. Both connecting
homomorphisms are given by projections of the restrictions of the
convolution product. In light of this observation, the $H$-cleft
property of the extension $B\subseteq A$ is equivalent to the
existence of elements $j$ and $\tildej$ in the two bimodules of
the Morita context, such that the connecting homomorphisms map
$j\otimes \tildej$ and $\tildej\otimes j$, respectively, to
the unit element of the convolution algebra. Inspired by this
example, we term an $L$-$\cC$ bicomodule $\Sigma$, for a pure coring 
extension $(\cD:L)$ of $(\cC:A)$, to be {cleft} if there exist
elements $j$ and $\tildej$ in the two bimodules of the Morita
context, associated to $\Sigma$ in \seref{contextext}, such that
the connecting homomorphisms map $j\otimes \tildej$ and
$\tildej\otimes j$, respectively, to the unit element of the
appropriate algebra in the Morita context. In Section \ref{sec:ex}
we show that most examples of `cleft extensions' -- i.e.
cleft algebra extensions by Hopf algebras or pure
Hopf algebroids, cleft extensions
by partial group actions, as well as cleft (weak) entwining structures (over
arbitrary base) -- are covered by this definition. 
It has to be noted that, in contrast to cleft extensions by Hopf
algebras (or just coalgebras), which determine cleft entwining structures,
cleft extensions (in the sense of \cite{BohmBrz:cleft}) by pure Hopf
algebroids correspond to pure coring extensions which do not arise from any
entwining structure. 
Cleft extensions by arbitrary (not necessarily pure) Hopf algebroids are not
even known to correspond to cleft bicomodules for any pure coring
extension. Instead, in Remark \ref{rem:non.pure} we construct a Morita
context of slightly different kind, such that the cleft property of an algebra
extension by any Hopf algebroid corresponds to the existence of mutually
inverse elements $j$ and $\tildej$ in its bimodules.

As it is well known, an extension $B\subseteq A$ of $k$-algebras is an
{\em $H$-cleft} extension for a Hopf algebra $H$ if and only if it is {\em
$H$-Galois} with {\em normal basis property}. The Galois condition
means that $A$ is a right $H$-comodule algebra with coinvariants
$B$, and $A\otimes_{B} A$ is isomorphic (in the category of left
$A$-modules, right $H$-comodules) to $A\otimes_k H$ via the so
called canonical isomorphism. The normal basis property means that
$A$ is isomorphic to $B\ot_k H$ in the category of left
$B$-modules, right $H$-comodules.

The notion of Galois extensions went through a wide generalization in
the past years. It was initiated by an observation by Brzezi\'nski
\cite{Brz:str} that if $B\subseteq A$ is an $H$-Galois extension then
$\cC\colon = A\ot_k H$ is a {\em Galois $A$-coring}. This means that $A$ is
a $\cC$-comodule with $\cC$-coinvariants $B$, and $\cC$ is isomorphic to the
canonical Sweedler's
coring $A\sstac{B} A$ via the canonical isomorphism.

As a next step of generalization \cite{ElK:comcor},
El Kaoutit and G\'omez-Torrecillas introduced {\em Galois comodules}
$\Sigma$ for an $A$-coring $\cC$ as right $\cC$-comodules that are
finitely generated and projective as $A$-modules and $\cC$ is
isomorphic to a comatrix coring
$\Hom_A(\Sigma,A)\sstac{T} \Sigma$ via the canonical isomorphism, where the
notation $T\colon = \End^\cC(\Sigma)$ is used. The existence of
such a Galois comodule $\Sigma$ ensures that the comonad functor
$-\sstac{A} \cC$ on the category $\cM_A$ of right $A$-modules comes
from an adjunction of functors 
\begin{equation}\label{eq:adjoints}
-\otimes_T \Sigma:\cM_T\to \cM_A\quad\textrm{and}\quad
{\rm Hom}_A(\Sigma,-):\cM_A\to \cM_T\ .
\end{equation}
Hence, by standard Eilenberg-Moore
type arguments (see e.g. \cite[VI.3, Theorem 1]{MacL}), the diagram
\beq\label{fig:EiMo}
\xymatrix
{
\cM_A \ar@<-.5ex>[dd]_{-\sstac{A}\cC}\ar@<.5ex>[ddrrr]^{{\rm
Hom}_A(\Sigma,-)}&&&\\
&\\
\cM^\cC\ar@<-.5ex>[uu]_{F^\cC}&&&\cM_T
\ar@<.5ex>[uulll]^{-\sstac{T} \Sigma}
\ar[lll]^{-\sstac{T} \Sigma}
}
\eeq
is a commutative diagram in the 2-category of categories, in the sense
that the inner triangle strictly commutes (here $\cM^\cC$ denotes the category
of right $\cC$-comodules and $F^\cC$ denotes
the forgetful functor) and the outer one does upto the natural
isomorphism
\beq\label{eq:can}
\can_N:\  \Hom_A(\Sigma,N)\stac{T} \Sigma \to N\stac{A} \cC,
\qquad
\phi_N\stac{T} x \mapsto \phi_N(x^{[0]})\stac{A} x^{[1]},
\eeq
for any right $A$-module $N$, where $x\mapsto x^{[0]}\stac{A} x^{[1]}$ denotes
the $\cC$-coaction in $\Sigma$.

Finally, in \cite{Wis:galcom} Wisbauer relaxed the finitely generated
projectivity 
condition and defined Galois comodules $\Sigma$ for an $A$-coring
$\cC$ with the requirement that the diagram \eqref{fig:EiMo} is
commutative in the above sense, i.e. such that \eqref{eq:can} is a natural
isomorphism. 
Throughout the paper Galois comodules are meant in this most general sense.

As one of our main results, it is shown in \coref{jJ}
that an $L$-$\cC$ bicomodule $\Sigma$, for a pure coring extension 
$(\cD,L)$ of $(\cC,A)$, is cleft if and only if it is a Galois
comodule for $\cC$ (in the sense of \cite{Wis:galcom}) and
isomorphic to $T \otimes_L \cD$ as a left
module for the algebra $T\colon ={\rm End}^\cC(\Sigma)$ and as a comodule
for the coring $\cD$ (i.e. the {\em normal basis} property holds).

Let $\Sigma$ be a Galois comodule for an $A$-coring $\cC$ in the
sense of \cite{Wis:galcom}, and $T\colon ={\rm End}^\cC(\Sigma)$.
Then the Eilenberg-Moore comonad for the adjunction \eqref{eq:adjoints}
is naturally equivalent to the comonad $-\sstac A\cC$, and the 
comparison functor is naturally equivalent to $-\otimes_T \Sigma 
:\cM_T\to \cM^\cC$ (cf. diagram \eqref{fig:EiMo}). Hence $-\sstac T
\Sigma:\cM_T\to \cM_A$ is {\em comonadic} (or {\em tripleable}) if and
only if the functor in the bottom row of diagram \eqref{fig:EiMo}
is an equivalence (i.e. the {\em Strong Structure Theorem} holds).
Note that this functor in the bottom row of diagram
\eqref{fig:EiMo} possesses a right adjoint (cf.
\cite[18.21]{BrzWis:cor}), the functor \beq \label{eq:coinvfun}
\Hom^\cC(\Sigma,-):\cM^\cC\to \cM_T. \eeq In \thref{surjectivity} (1)
and Theorem 
\ref{thm:w.str.thm} we show that if $\Sigma$ is an $L$-$\cC$
bicomodule for a pure coring extension $(\cD:L)$ of $(\cC:A)$, then the
surjectivity of one of the connecting homomorphisms in the
associated Morita context implies both the Galois property of
$\Sigma$ as a right $\cC$-comodule and the fully faithfulness of
the functor \eqref{eq:coinvfun} (i.e. the {\em Weak Structure
Theorem}). The Strong Structure Theorem does not follow even by
the strictness of the Morita context without further assumptions.
In Theorem \ref{thm:s.str.thm} we find sufficient conditions for
it to hold. Since Morita contexts associated to cleft bicomodules
are strict, these theorems -- extending
\cite[Proposition 4.8]{Abuh:cleftentw}, \cite[Theorem
4.5]{CaVerWang:cleftentw} and \cite[Theorems 4.9 and 4.10]{Abuh:cleftentw} --
hold for them. 
Using similar ideas, we prove a Strong Structure Theorem also for cleft
extensions by arbitrary Hopf algebroids (\thref{thm:nonpure.str.str.thm}). It
yields in particular a {\em corrected form} of \cite[Theorem
  4.2]{Bohm:hgdint}, whose proof in the original journal version turned out to
be incorrect.

\bigskip

{\bf Notational conventions.} All algebras of the paper are
associative unital algebras over a fixed commutative ring $k$. The
multiplication in an algebra $R$ will be denoted by $\mu_R$ and
the unit element by $1_R$. 
The algebra with the same $k$-module structure and opposite multiplication is
denoted by $R^{op}$.
For the category of right (resp. left)
$R$-modules we use the symbol $\cM_R$ (resp. ${}_R \cM$) and for
its hom sets we write $\Hom_R(-,-)$ (resp. ${}_R\Hom(-,-)$). The
forgetful functor $\cM_R\to \cM_k$ will be denoted by $G^R$.

$R$-rings over an algebra $R$ are monoids in the monoidal category 
of $R$-$R$ bimodules. An $R$-ring is determined by a pair, consisting 
of an algebra $A$ and an algebra homomorphism (the unit of the monoid) 
$\eta_A:R\to A$.

Corings (co-rings) over an algebra $R$ are meant to be comonoids
in the monoidal category of $R$-$R$ bimodules. That is, an
$R$-coring $\cC$ consists of an $R$-$R$ bimodule (denoted by the
same symbol $\cC$), an $R$-$R$ bilinear coassociative coproduct
$\Delta_\cC$ and an $R$-$R$ bilinear counit $\epsilon_\cC$ (see
\cite[17.1]{BrzWis:cor}). 
For the coproduct of an element $c\in \cC$ we use Sweedler-Heyneman index
notation, i.e. write $\Delta_\cC(c)=c_{(1)}\sstac R c_{(2)}$, without 
denoting summation explicitly.
The category of right (resp. left)
$\cC$-comodules (cf. \cite[18.1-2]{BrzWis:cor}) will be denoted by
$\cM^\cC$ (resp. ${}^\cC\cM$) and its hom sets by $\Hom^\cC(-,-)$
(resp. ${}^\cC\Hom(-,-)$). The forgetful functor $\cM^\cC\to
\cM_R$ will be denoted by $F^\cC$. Composition of morphisms in the 
$k$-linear category $\cM^\cC$ equips the set $\End^\cC(\Sigma)$ of
endomorphisms of an object $\Sigma$ with an algebra structure. Symmetrically,
also the composition of morphisms in ${}^\cC\cM$ makes the 
set ${}^\cC\End(\Lambda)$ of endomorphisms of $\Lambda\in
{}^\cC\cM$ an algebra.

\section{Morita contexts associated to comodules}\label{sec:Morita}
Generalizing constructions in \cite{ChaSwe:hagal}, \cite{Doi:smmor},
\cite{Abuh:cleftentw} and \cite{CaVerWang:cleftentw}, 
Caenepeel, De Groot and Vercruysse in \cite[Section 4]{CaDeGrVe:comcor} 
associated Morita contexts to comodules of an $A$-coring $\cC$.
For any right $\Cc$-comodule $\Sigma$, they constructed a Morita context, 
connecting the algebras $\End^\Cc(\Sigma)$ and the right dual $\Cc^*$ of
$\cC$ (cf. ${\mathbb M}'(\Sigma)$ below). Dually, for a left $\Cc$-comodule
$\Lambda$ they constructed a Morita context, connecting
${^\Cc\End}(\Lambda)^{\rm op}$ and the left dual, $\*C$ (cf. ${\mathbb
M}'(\Lambda)$ below). In the case of a right $\Cc$-comodule $\Sigma$,
which is finitely generated and projective as right $A$-module, these
constructions yield a Morita context connecting $\End^\Cc(\Sigma)$ and $\*C$
(cf. ${\mathbb M}'(\Sigma^*)$ below). 
In the present section we give a generalization of this last Morita context to
arbitrary right $\cC$-comodules $\Sigma$. 
For generalities about Morita contexts, we refer to \cite[Chapter II.3]{Bass}, 
where they are termed sets of (pre-)equivalence data.

Recall (e.g. from \cite[17.8 and 19.1]{BrzWis:cor}, where however the
  convention of opposite multiplication is used,) 
that for an $A$-coring 
$\cC$ the $k$-module ${}^*\cC\colon ={}_A\Hom(\cC,A)$ is an $A$-ring with
multiplication
$$
(ff')(c)=f'\big(c^{(1)} f(c^{(2)})\big),\qquad
\textrm{for } f,f'\in {}^*\cC,\ c\in \cC,
$$
and unit map $A\to {}^*\cC$, $a\mapsto \big(\ c\mapsto \epsilon_\cC(ca)\
\big)$. 
Any right $\cC$-comodule $\Sigma$ possesses a canonical right ${}^*\cC$-module
structure with action
$$x f=x^{[0]} f(x^{[1]}),\qquad
\textrm{for } f\in {}^*\cC,\ x\in \Sigma.
$$
$\cC$-colinear maps are ${}^*\cC$-linear. In a symmetric way, also the right
dual $\cC^*=\mathrm{Hom}_A(\cC,A)$ has an $A$-ring structure and there is a
faithful functor ${}^\cC\cM\to {}_{\cC^*}\cM$.

Recall from \cite[Section 4]{CaDeGrVe:comcor} that, for a right comodule
$\Sigma$ of an $A$-coring $\cC$ and $T\colon =\mathrm{End}^\cC(\Sigma)$, 
$\mathrm{Hom}^\cC(\cC,\Sigma)$ is a $T$-$\cC^*$ bimodule with
$$
(t\omega g)(c)=t\circ\omega (g( c^{(1)}) c^{(2)}),
\quad \qquad \textrm{for } t\in T, \ \omega \in \mathrm{Hom}^\cC(\cC,\Sigma),\
g\in \cC^*, 
$$
and $\Sigma^*\colon=\mathrm{Hom}_A(\Sigma,A)$ is a $\cC^*$-$T$ bimodule with 
$$
(g \xi t)(x)=g( \xi\circ t(x^{[0]})x^{[1]}),
\quad\qquad\textrm{for } t\in T, \ \xi \in \Sigma^*,\ g\in
\cC^*.\qquad\qquad 
$$
They constitute a Morita context
\begin{equation}\label{eq:OldMor}
{\mathbb M}'(\Sigma)=
\big(\ T,\cC^*,
\mathrm{Hom}^\cC(\cC,\Sigma),\Sigma^*,\newblacktriangle',\newtriangle'\ \big), 
\end{equation}
with connecting homomorphisms
\begin{eqnarray*}
&\newblacktriangle':\Sigma^*\stac{T} \mathrm{Hom}^\cC(\cC,\Sigma)\to
\cC^*,\qquad  
&\xi\stac{T} \omega\mapsto \xi\circ \omega,\\
&\newtriangle':\mathrm{Hom}^{\cC}(\cC,\Sigma)\stac{\cC^*} \Sigma^*\to T,\qquad
&\omega\stac{\cC^*}\xi \mapsto 
\big(\ x\mapsto \omega(\xi(x^{[0]})x^{[1]})\ \big).
\end{eqnarray*}
For the purposes of the present paper, however, another Morita context
associated to $\Sigma\in \cM^\cC$ turns
out to be more useful. Define the $k$-module 
\beq\label{eq:Q}
Q\colon =\{\ q\in \Hom_A(\Sigma,{}^*\cC)\ |
\ \forall x\in \Sigma,c\in \cC\quad
q(x^{[0]})(c)x^{[1]}=c^{(1)}q(x)(c^{(2)})\ \}.
\eeq
In the following lemma some properties of the $k$-module $Q$ are collected,
needed in order to see that $\Sigma$ and $Q$ are bimodules connecting the
algebras $T$ and ${}^*\cC$.
\begin{lemma}\label{lem:Qprop}\lelabel{Qprop}\lelabel{3.4}
For an $A$-coring $\cC$ and its right comodule $\Sigma$, the $k$-module $Q$ in
\eqref{eq:Q} obeys the following properties.
\begin{zlist}
\item $Q$ is isomorphic to the $k$-module 
\begin{equation}\eqlabel{defQ'}
Q'\colon=\{\ q\in {}_A{\rm Hom}(\cC,\Sigma^*)\ |\
\forall x\in \Sigma, c\in \cC\quad 
c^{(1)}q(c^{(2)})(x)=q(c)(x^{^{[0]}})x^{[1]}\ \},
\end{equation}
defined in terms of the left $A$-module  $\Sigma^*\colon ={\rm
  Hom}_A(\Sigma,A)$;
\item Let $M$ be a right $\cC$-comodule.
For any $q\in Q$ and $m\in M$, the map
$\Sigma\to M$, $x\mapsto mq(x)$
is right $\cC$-colinear, i.e. for every $m\in M$
there is a $k$-linear map 
$$
Q\to \Hom^\Cc(\Sigma,M),\quad q\mapsto m q(-);
$$
\item $Q$ is a $k$-submodule of $\Hom_{{}^*\cC}(\Sigma,{}^*\cC)$;
\item $Q$ is a ${}^*\cC$-$T$ bimodule, for $T\colon =\End^\cC(\Sigma)$, with
  actions  
\begin{eqnarray*}
(fq)(x)&\colon = f q(x),\qquad &\textrm{for } f\in {}^*\cC,\ q\in Q,\ x\in
  \Sigma\qquad \textrm{and}\\
(qt)(x)&\colon =q\big(t(x)\big),\qquad &\textrm{for }q\in Q,\ \ t\in T,\,\
  x\in \Sigma. 
\end{eqnarray*}
\end{zlist}
\end{lemma}
\begin{proof}
\ul{(1)} The isomorphism is given by switching the arguments. That is, by the
map 
$$
Q\to Q',\qquad q\mapsto \big(\ c\mapsto q(-)(c)\ \big).
$$

{\underline{(2)}}
Since ${}^*\cC$ is an $A$-ring and the elements of $Q$ are right $A$-linear,
the map $\Sigma\to M$, $x\mapsto mq(x)$
is right $A$-linear,
for $q\in Q$ and a right $\cC$-comodule $M$ and $m\in M$. In order to see that
it is also right $\cC$-colinear, 
use the right $A$-linearity of a $\cC$-coaction in the
first equality and \eqref{eq:Q} in the second one, to conclude that,
for any right $\cC$-comodule $M$, $m\in M$, $q\in Q$ and $x\in \Sigma$, 
\begin{eqnarray*}
\big(m^{[0]} q(x)(m^{[1]})\big)^{[0]}\stac{A}(m^{[0]}
q(x)(m^{[1]})\big)^{[1]} 
&=&m^{[0]}\stac{A}m^{[1]} q(x)(m^{[2]})\\
=m^{[0]}\stac{A} q(x^{[0]})(m^{[1]})x^{[1]}
&=& m^{[0]}q(x^{[0]})(m^{[1]})\stac{A}x^{[1]}.
\end{eqnarray*}

{\underline{(3)}}
Using the right $A$-linearity of $q\in Q$ and the defining identity
\eqref{eq:Q}, one checks that, for $x\in\Sigma$, $f\in {}^*\cC$ and $c\in
\cC$, 
\begin{eqnarray*}
q(xf)(c)
&=&q\big(x^{[0]}f(x^{[1]})\big)(c)
=q(x^{[0]})(c)f(x^{[1]})=f\big(q(x^{[0]})(c)x^{[1]}\big) \\
&=&f\big(c^{(1)} q(x)(c^{(2)})\big)=\big(q(x)f\big)(c),
\end{eqnarray*}
where in the first equality the form of the ${}^*\cC$-action in $\Sigma$ has
been used, and in the last one the multiplication law in ${}^*\cC$.

{\underline{(4)}}
For $f\in {}^*\cC$ and $q\in Q$, $fq$ is an element of 
$\Hom_A(\Sigma,{}^*\cC)$ by the right $A$-linearity of $q$ and the fact that
${}^*\cC$ is an $A$-ring. Since $q$ is an element of $Q$, so is $fq$ as, for
$x\in\Sigma$ and $c\in \cC$, 
\begin{eqnarray*}
(fq)(x^{[0]})(c)x^{[1]}
&=&\big(fq(x^{[0]})\big)(c)x^{[1]}
=q(x^{[0]})\big(c^{(1)} f(c^{(2)})\big)x^{[1]}   \\
&=&c^{(1)} q(x)\big(c^{(2)} f(c^{(3)})\big)
=c^{(1)}\big(fq(x)\big)(c^{(2)})
=c^{(1)} (fq)(x)(c^{(2)}),
\end{eqnarray*}
where the first and last equalities follow by the the form of the
${}^*\cC$-action in $Q$ and the second and penultimate equalities follow by
the multiplication law in ${}^*\cC$. The third equality follows by
\eqref{eq:Q}. 

Since $q\in Q$ and $t\in T$ are right $A$-linear, $qt$ is an element of
$\Hom_A(\Sigma,{}^*\cC)$. Since $t\in T$ is colinear and $q$ is an element of
$Q$, 
it follows that $qt\in Q$ as, for $x\in\Sigma$ and $c\in \cC$, 
$$
(qt)(x^{[0]})(c)x^{[1]}=q\big(t(x)^{[0]}\big)(c)t(x)^{[1]}
=c^{(1)} (qt)(x)(c^{(2)}),
$$
where the form of the $T$-action in $Q$ has been used.

It is straightforward to check that both actions are associative and unital
and they commute.
\end{proof}

\begin{remark}\label{rem:Qfgp}
In the case when $\cC$ is a locally projective left $A$-module, the 
$k$-module \eqref{eq:Q} has a particularly simple characterization, as
$Q\equiv \mathrm{Hom}_{{}^*\cC}(\Sigma,\*C)$. Indeed, by Lemma \ref{lem:Qprop} 
(3),  
$Q\subseteq \mathrm{Hom}_{{}^*\cC}(\Sigma,\*C)$. The converse inclusion is
proven as follows.
Recall that local projectivity of
the left $A$-module $\cC$ means that for any finite subset ${\mathcal
  S}\subset \Cc$, there exists a dual basis $\{e_i\}\subset \cC$ and
$\{f_i\}\subset \*C$ such that $c=\sum_if_i(c)e_i,$ for all 
$c\in{\mathcal S}$. For an element $x\in \Sigma$, fix finite sets
$\{x_j\}\subset \Sigma$ and $\{c_j\}\subset \cC$ such that
$$
x^{[0]}\stac A x^{[1]}=\sum_{j=1}^s x_j\stac A c_j.
$$
Take $q\in \mathrm{Hom}_{{}^*\cC}(\Sigma,\*C)$ and $c\in \cC$ and introduce
the set  
${\mathcal S}\colon= \{c_1,\cdots,c_s, c^{(1)}q(x)(c^{(2)})\}$ $\subset \Cc$. 
By the assumption of local projectivity, there exists a dual basis 
$\{ e_i\}\subset \cC$ and $\{f_i\}\subset {}^*\cC$ associated to ${\mathcal
  S}$, hence we can write 
\begin{eqnarray*}
q(x^{[0]})(c)x^{[1]}
&=&\sum_i q(x^{[0]})(c)f_i(x^{[1]})e_i
= \sum_i q(x^{[0]}f_i(x^{[1]}))(c)e_i\\
&=&\sum_i q(x f_i)(c)e_i
=\sum_i\big(q(x)f_i\big)(c)e_i\\
&=&\sum_i f_i\big(c^{(1)}q(x)(c^{(2)})\big)e_i=
c^{(1)}q(x)(c^{(2)}),
\end{eqnarray*}
where the fourth equality follows by the right ${}^*\cC$-linearity of $q$.
This shows that $q$ belongs to the $k$-module $Q$ in \eqref{eq:Q}.

If the $A$-coring $\cC$ is locally projective as a left $A$-module then
the image of a right $\Cc$-comodule $\Sigma$ under an element $q\in
\Hom_{{}^*\cC}(\Sigma,{}^*\cC)$ lies within the rational part $({}^*\cC)^{\rm
rat}$ of ${}^*\cC$, cf. \cite[20.1]{BrzWis:cor}. (Recall that the rational
part of a right ${}^*\cC$-module is the biggest ${}^*\cC$-submodule
that possesses a $\Cc$-comodule structure.) This implies that
$Q=\Hom_{^*\cC}(\Sigma,{}^*\cC)=\Hom_{^*\cC}(\Sigma,({}^*\cC)^{\rm
rat})=\Hom^\Cc(\Sigma,({}^*\cC)^{\rm rat})$. 

If $\cC$ is a finitely generated and projective left $A$-module, with
dual bases $\{c_i\}\subset \cC$ and $\{f_i\}\subset {}^*\cC$, then
${}^*\cC$ possesses a right $\cC$-comodule structure with coaction
$f\mapsto \sum_i f f_i\sstac{A} c_i$. In this case $Q$ is identical to
the $k$-module $\Hom^\cC(\Sigma,{}^*\cC)\equiv
 \Hom_{{}^*\cC}(\Sigma, {}^*\cC)$. 
\end{remark}
In light of \leref{Qprop}, there is another Morita context associated to
$\Sigma$,  
\begin{equation}\label{eq:SiMor}
{\mathbb M}(\Sigma)=(T,{}^*\cC,\Sigma,Q,\newblacktriangle,\newtriangle),
\end{equation}
where $txf=t\big(x^{[0]}f(x^{[1]})\big)$, for $t\in T$, $x\in \Sigma$ and $f\in
{}^*\cC$, with connecting homomorphisms 
\begin{eqnarray}
\newblacktriangle:Q\stac{T} \Sigma &\to {}^*\cC,\qquad q\stac{T} x &\mapsto
q(x),\label{eq:F}\eqlabel{defF}\\
\newtriangle:\Sigma\stac{{}^*\cC} Q&\to T,\qquad x\stac{{}^*\cC} q&\mapsto
xq(-).\label{eq:G}\eqlabel{defG}
\end{eqnarray}
In a symmetric way, to a left $\cC$-comodule $\Lambda$ one can associate
Morita contexts ${\mathbb M}'(\Lambda)$, connecting the algebras
${}^\cC\mathrm{End}(\Lambda)^{op}$ and ${}^*\cC$, and ${\mathbb
  M}(\Lambda)$, connecting ${}^\cC\mathrm{End}(\Lambda)^{op}$ and $\cC^*$.
\begin{remark}\label{rem:Sifgp}
Note that a finitely generated and projective right $A$-module $\Sigma$ is a
right comodule for an $A$-coring $\cC$ if and only if $\Sigma^*$ is a left
$\cC$-comodule. Indeed, with the help of a dual basis $\{x_i\}\subset \Sigma$
and $\{\xi^i\}\subset \Sigma^*$, a bijective correspondence is given between
right $\cC$-coactions $x\mapsto x^{[0]}\sstac{A} x^{[1]}$ in $\Sigma$, and 
left $\cC$-coactions $\xi\mapsto \sum_i \xi(x_i^{[0]})x_i^{[1]}\sstac{A}
\xi^i$ in $\Sigma^*$. In this case the
Morita context ${\mathbb M}(\Sigma)$ coincides with 
${\mathbb M}'(\Sigma^*)$ and  ${\mathbb M}'(\Sigma)$ coincides with 
${\mathbb M}(\Sigma^*)$.
\end{remark}
Since any right $\cC$-comodule $\Sigma$ has also a right $\*C$-module
structure, one can associate a further Morita context with it, as in
\cite[II.4]{Bass}. Namely,  
\begin{equation}\eqlabel{tildeQ}
{\mathbb N}(\Sigma)=
\big(\End_{\*C}(\Sigma),\*C,\Sigma,\Hom_{\*C}(\Sigma,\*C),\newblacktriangleup,
\newtriangleup \big) 
\end{equation}
with connecting maps 
\begin{eqnarray}
\newblacktriangleup:&
\Hom_{\*C}(\Sigma,\*C)\stac{\End_{\*C}(\Cc)}\Sigma \to \*C,\qquad
&q\stac{\End_{\*C}(\Cc)}x\mapsto q(x)\eqlabel{deftildeF}\\ 
\newtriangleup:&
\Sigma\stac{\*C}\Hom_{\*C}(\Sigma,\*C)\to \End_{\*C}(\Cc),\qquad
&x\ \stac{\*C}\ q\ \mapsto \ xq(-).\eqlabel{deftildeG} 
\end{eqnarray}
In the next generalization of \cite[Proposition 4.7]{CaDeGrVe:comcor} the
relationship between the Morita contexts 
${\mathbb N}(\Sigma)$ in \equref{tildeQ} and ${\mathbb M}(\Sigma)$ in
\eqref{eq:SiMor} is investigated.

\begin{proposition}\prlabel{2.6}
Let $\Cc$ be an $A$-coring and $\Sigma$ a right $\Cc$-comodule. There exists a
morphism of Morita contexts ${\mathbb M}(\Sigma)\to {\mathbb N}(\Sigma)$,
which becomes an isomorphism if $\Cc$ is locally projective as left $A$-module.
\end{proposition}

\begin{proof}
There exist inclusions $T=\End^\Cc(\Sigma)\subset \End_{\*C}(\Sigma)$ and, 
by \leref{Qprop} (3), $Q \subset\Hom_{\*C}(\Sigma,\*C)$.
Comparing \equref{defF} and \equref{defG} with \equref{deftildeF} and
\equref{deftildeG}, it is straightforward to see that these inclusions,
together with 
the identity maps of $\*C$ and $\Sigma$, establish a morphism of Morita
contexts. 

Now assume that $\Cc$ is locally projective as left $A$-module. By
\cite[19.2-3]{BrzWis:cor}, for any right $\cC$-comodules $M$ and $M'$,
$\mathrm{Hom}^\cC(M,M')=\mathrm{Hom}_{{}^*\cC}(M,M')$. Hence in 
particular $T=\mathrm{End}_{{}^*\cC}(\Sigma)$. By Remark \ref{rem:Qfgp},
$Q=\mathrm{Hom}_{{}^*\cC}(\Sigma,\*C)$, which completes the proof.
\end{proof}

By standard Morita theory, if the connecting map $\newblacktriangle$ in the
Morita context ${\mathbb M}(\Sigma)$
in \eqref{eq:SiMor} is surjective, 
then $\Sigma$ is a finitely generated projective left $T$-module and a right
${}^*\cC$-generator. If $\newtriangle$ is surjective 
then $\Sigma$ is a finitely generated projective right
${}^*\cC$-module and a left $T$-generator. 
\begin{lemma}\label{lem:Afgp}
Let $\cC$ be an $A$-coring and $\Sigma$ a right comodule. Consider the Morita 
context ${\mathbb
  M}(\Sigma)=(T,{}^*\cC,\Sigma,Q,\newblacktriangle,\newtriangle)$ 
in \eqref{eq:SiMor}.   
\begin{zlist}
\item If the connecting map $\newblacktriangle$ in \equref{defF} is surjective
  then $\cC$ is a finitely generated projective left $A$-module.
\item If the connecting map $\newtriangle$ in \equref{defG} is surjective then
  $\Sigma$ is a 
finitely generated projective right $A$-module.
\end{zlist}
\end{lemma}
\begin{proof}
\underline{(1)}
If $\newblacktriangle$ is surjective then there exist finite sets
$\{q_i\}\subset Q$ and $\{x_i\}\subset \Sigma$ such that 
$$
\epsilon_\cC=\sum_i q_i\newblacktriangled x_i\equiv 
\sum_i q_i(x_i).
$$
Introduce finite sets $\{f_j\}\subset {}^*\cC$ and $\{c_j\}\subset \cC$ via
the requirement that 
$$
\sum_j f_j\stac{A} c_j =\sum_i q_i(x_i^{[0]})\stac{A} x_i^{[1]}.
$$
We claim that they are dual bases. Indeed, for any $c\in \cC$,
$$
\sum_j f_j(c)c_j=\sum_i q_i(x_i^{[0]})(c) x_i^{[1]}=
\sum_i c^{(1)}q_i(x_i)(c^{(2)})=c^{(1)}\epsilon_\cC(c^{(2)})
=c,
$$
where the second equality follows by the definition \eqref{eq:Q} of $Q$.
Hence $\cC$ is a finitely generated projective left $A$-module, as
stated. 

\underline{(2)}
If $\newtriangle$ is surjective then there exist finite sets $\{x_i\}\subset
\Sigma$ and $\{q_i\}\subset Q$ such that 
$$
1_T=\sum_i x_i\newtriangled q_i\equiv \sum_i x_i^{[0]} q_i(-)(x_i^{[1]}).
$$
A dual basis for the right $A$-module $\Sigma$ can be constructed 
introducing finite sets $\{y_j\}\subset \Sigma$ and $\{\xi_j\}\subset
\Sigma^*$ via the requirement that
$$
\sum_j y_j\stac{A}\xi_j =\sum_i x_i^{[0]}\stac{A} q_i(-)(x_i^{[1]}).
$$
\end{proof}

\begin{theorem}\thlabel{tensorff}
Let $\cC$ be an $A$-coring and $\Sigma$ its right comodule. Consider the Morita
context ${\mathbb M}(\Sigma)$, associated to $\Sigma$ in
\eqref{eq:SiMor}. If the connecting map $\newtriangle$ in \equref{defG}
is surjective, then the functor $-\sstac T\Sigma:\cM_T\to \cM^\cC$ in the
bottom row of diagram \eqref{fig:EiMo} is fully faithful. 
\end{theorem}
\begin{proof}
The proof consists of a verification of the
bijectivity of the unit of the adjunction of functors $-\sstac T
\Sigma:\cM_T\to \cM^\cC$ and $\Hom^\cC(\Sigma,-):\cM^\cC\to \cM_T$,
i.e. the map 
\beq \label{eq:uN}
\eta_N\ :\ N\to \Hom^\cC(\Sigma, N\stac{T} \Sigma),\qquad
n\mapsto (\ x\mapsto n\stac{T} x\ ),
\eeq
for any right $T$-module $N$.
Choose elements $\{x_i\}\subset \Sigma$ and $\{q_i\}\subset Q$ such that
$\sum_i x_i\newtriangled q_i=1_T$. The inverse of the map \eqref{eq:uN}
can be constructed as 
$$
{\widetilde \eta}_N\ :\ \Hom^\cC(\Sigma, N\stac{T} \Sigma)\to N, \qquad
\zeta_N \mapsto (N\stac{T} \newtriangle)\big(\sum_i \zeta_N(x_i)\stac{{}^*\cC}
q_i\big). 
$$
Indeed, the identity ${\widetilde \eta}_N\circ \eta_N=N$ obviously holds
true. For the other equality, use the associativity of the Morita context
${\mathbb M}(\Sigma)$ to compute, for $\zeta_N\in
\mathrm{Hom}^\cC(\Sigma,N\sstac{T} \Sigma)$ and $x\in \Sigma$,
\begin{eqnarray*} 
\big(\eta_N \circ {\widetilde \eta}_N\big)(\zeta_N)(x)
&=&(N\otimes_T\newtriangle\sstac{T} \Sigma)
(\sum_i \zeta_N(x_i)\otimes_{\*C}q_i \otimes_T x)\\
&=&(N\otimes_T\Sigma \sstac{{}^*\cC} \newblacktriangle)
(\sum_i \zeta_N(x_i)\otimes_{\*C}q_i \otimes_T x)\\
&=&\sum_i \zeta_N(x_i)(q_i\newblacktriangled x)
=\zeta_N(\sum_ix_i(q_i\newblacktriangled x))
=\zeta_N(\sum_i(x_i\newtriangled q_i)x)=\zeta_N(x),
\end{eqnarray*}
where in the fourth equality we used the
right ${}^*\cC$-linearity of $\zeta_N\in \Hom^\cC(\Sigma, N\stac{T} \Sigma)$.
\end{proof}

\begin{proposition}\label{prop:Cfgp}
Let $\cC$ be an $A$-coring which is finitely generated and projective as a left
$A$-module and let $\Sigma$ be a right comodule. Put  $T\colon
=\End^\cC(\Sigma)$.
The following
statements are equivalent.
\begin{zlist}
\item The Morita context ${\mathbb M}(\Sigma)$, associated to $\Sigma$ in
  \eqref{eq:SiMor}, is strict;  
\item $\Sigma$ is a Galois comodule and finitely generated
and projective 
as a right $A$-module and faithfully flat as left $T$-module;
\item The functor $-\sstac{T}
\Sigma:\cM_T\to \cM^\cC$, in the bottom row of diagram
\eqref{fig:EiMo}, is an equivalence with inverse \eqref{eq:coinvfun}.
\end{zlist}
\end{proposition}
\begin{proof}
\underline{(1)$\Rightarrow$ (3)}
If the Morita context is strict then the functors
$-\sstac{T}\Sigma: \cM_T\to \cM_{{}^*\cC}$ and $-\sstac{{}^*\cC}
Q:\cM_{{}^*\cC}\to \cM_T$ are inverse equivalences. 
Since both functors $-\sstac{{}^*\cC} Q$ and 
$\Hom_{{}^*\cC}(\Sigma,-):\cM_{{}^*\cC}\to \cM_T$ are right
adjoints to $-\sstac{T}\Sigma$, by uniqueness of adjoint functors upto
natural equivalence, both are inverses of it. 
Since $\cC$ is a finitely generated and projective left
$A$-module, the categories $\cM^\cC$ and $\cM_{{}^*\cC}$ are
isomorphic (cf. \cite[19.6]{BrzWis:cor}).
Consequently, $-\sstac{T}\Sigma: \cM_T\to \cM^\cC$ is an equivalence
with inverse \eqref{eq:coinvfun}.

\underline{(3)$\Rightarrow$ (1)}
Since $\cC$ is a finitely generated and projective left
$A$-module, the categories $\cM^\cC$ and $\cM_{{}^*\cC}$ are
isomorphic, hence $-\sstac{T}\Sigma: \cM_T\to \cM_{{}^*\cC}$ and
$\Hom_{{}^*\cC}(\Sigma,-):\cM_{{}^*\cC}\to \cM_T$ are inverse equivalences.
The canonical strict Morita context, associated to
this equivalence, is equal to ${\mathbb M}(\Sigma)$ in \eqref{eq:SiMor}.

\underline{(2)$\Rightarrow$ (3)}
This follows by \cite[18.27 (2) (a)$\Rightarrow$ (c)]{BrzWis:cor}. 

\underline{(3)$\Rightarrow$ (2)}
Since we have proven already that (3) implies (1), 
$\Sigma$ is a finitely generated and projective 
right $A$ module by Lemma \ref{lem:Afgp} (2). Then it 
is a Galois comodule and a faithfully flat left $T$-module by \cite[18.27 (2)
  (c)$\Rightarrow$ (a)]{BrzWis:cor}.  
\end{proof}

\section{Morita theory for coring extensions}\selabel{contextext}
Let $\cD$ be a coring over the base $k$-algebra $L$ and $\cC$ a
coring over the $k$-algebra $A$. Assume that $\cC$ is a
$\cC$-$\cD$ bicomodule with the left regular $\cC$-coaction
$\Delta_\cC$ and some right $\cD$-coaction $\tau_\cC$. By
definition \cite[22.1]{BrzWis:cor}, this means that $\tau_\cC$ is
left $A$-linear (hence $\cC\otimes_{A} \cC$ is also a
$\cD$-comodule with coaction $\cC\otimes_{A}\tau_\cC$) and the
coproduct $\Delta_\cC$ is right $\cD$-colinear. Equivalently, the
coproduct $\Delta_\cC$ is right $L$-linear (hence
$\cC\otimes_{L}\cD$ is a left $\cC$-comodule with coaction
$\Delta_\cC\otimes_{L}\cD$) and the $\cD$-coaction $\tau_\cC$ is
left $\cC$-colinear. This situation was termed by Brzezi\'nski 
in \cite[Definition 2.1]{Brz:corext} as $\cD$ is a {\em right extension} 
of $\cC$.
A right coring extension $\cD$ of $\cC$ is said to be {\em pure} if the
equaliser
\begin{equation}\label{eq:M}
\xymatrix{
M \ar[rr]^-{\varrho}&&
M \ot_A {\mathcal C} \ar@<2pt>[rr]^-{\varrho\ot_A {\mathcal C}}
\ar@<-2pt>[rr]_-{M\ot_A \Delta_{\mathcal C}}&&
M \ot_A {\mathcal C} \ot_A {\mathcal C}
}
\end{equation}
in $\cM_L$ is $\cD\ot_L \cD$-pure, i.e. it is preserved by the functor $-
\ot_L \cD\ot_L \cD:\cM_L\to \cM_L$, for any right $\cC$-comodule $(M,\varrho)$.

Based on \cite[22.3]{BrzWis:cor} and its Erratum, 
the following can be proven
(see also the arXiv version of \cite[Theorem 2.6(2)]{Brz:corext}).
If ${\mathcal D}$ is a pure right coring extension of ${\mathcal C}$
then there is a $k$-linear functor $U=  - \Box_{\mathcal C} \, {\mathcal
  C}:{\mathcal M}^{\mathcal C} \to {\mathcal M}^{\mathcal D}$ that makes
commutative the diagram 
$$
\xymatrix
{\cM^\cC\ar[rrrr]^{U}\ar[rd]^{F^\cC}&&&&\cM^\cD\ar[ld]_{F^\cD}\\
&\cM_A\ar[rd]^{G^A} && \cM_L\ar[ld]_{G^L}&\\
&&\cM_k&&\ ,
}
$$
involving $U$ and four forgetful functors.
The functor $U$ has the following explicit form.
Using the right $\cD$-
coaction $\tau_\cC:c\mapsto c_{[0]}\otimes_{L}c_{[1]}$, for $c\in
\cC$, (note our convention to use character $\tau$ for
$\cD$-coactions and lower indices of the Sweedler type to denote
the components of the coproduct and of the coactions of the coring
$\cD$) 
the right $L$-action on $M$ in \eqref{eq:M} comes out as
$$
ml\colon = m^{[0]}\epsilon_\cC({m^{[1]}}l),
\qquad \textrm{for } m\in M \textrm{ and } l\in L.
$$
The cotensor product $M\Box_{\cC} \,\cC\cong M$ yields the
$\cD$-coaction
$$
\tau_M\ :\ M\to M\stac{L} \cD, \qquad m\mapsto m_{[0]}\stac{L} m_{[1]}\colon =
m^{[0]}\epsilon_\cC({m^{[1]}}_{[0]})\stac{L} {m^{[1]}}_{[1]},
$$
where $\varrho^M:m\mapsto m^{[0]}\otimes_{A} m^{[1]}$ denotes the
$\cC$-coaction in $M$ (note our convention to use character
$\varrho$ for $\cC$-coactions and upper indices of the Sweedler
type to denote the components of the coproduct and of the
coactions of the coring $\cC$). 
Note that $\tau_M$ is well defined and coassociative by the purity
assumption.
It is straightforward to check
that with this definition any right $\cC$-comodule map is
$\cD$-colinear. In particular, a right $\cC$-coaction, being
$\cC$-colinear by coassociativity, is $\cD$-colinear.

If $\Sigma$ is an object in the category ${}_L\cM^\cC$ of
$L$-$\cC$ bicomodules, i.e. it is a left $L$-module and a right
$\cC$-comodule with left $L$-linear $\cC$-coaction, then 
$T\colon=\mathrm{End}^\cC(\Sigma)$ is an $L$-ring with unit map $L\to T$,
$l\mapsto (\ x\mapsto lx\ )$.  
Furthermore, ${\rm Hom}^\cC(\Sigma,M)$ is a right $L$-module for any right
$\cC$-comodule $M$, via $(\varphi_Ml)(x)\colon = \varphi_M(lx)$, for
$\varphi_M\in {\rm Hom}^\cC(\Sigma,M)$, $l\in L$ and $x\in \Sigma$.
Hence, in addition to Brzezi\'nski's functor $U:\cM^\cC\to
\cM^\cD$, we can define another $k$-linear functor, 
$$V\colon = {\rm Hom}^\cC(\Sigma,-)\stac L \cD:\cM^\cC\to \cM^\cD.$$ 
Consider
the opposite of the category of $k$-linear functors and their
natural transformations. The full subcategory defined by the two
objects $U$ and $V$ determines a Morita context
\begin{equation}\label{eq:catMor} 
({\rm Nat}(V,V)^{op}\ ,\ {\rm Nat}(U,U)^{op}\ ,\ {\rm Nat}(V,U)\
,\ {\rm Nat}(U,V)\ ,\ \blacksquar\ ,\ \squar),
\end{equation}
where all algebra and bimodule structures are given by the
opposite composition of natural transformations and both
connecting homomorphisms $\blacksquar$ and $\squar$ are given by
projections of the restrictions of the opposite composition of
natural transformations.

In the following proposition an equivalent description of the Morita
context \eqref{eq:catMor}, in terms of sets of (co)module maps, is given.
In order to formulate it, we introduce a $k$-submodule of $Q'$, (i.e. the
$k$-module associated to
an $L$-$\cC$ bicomodule $\Sigma$ via
\equref{defQ'}). It is defined in terms of the $A$-$L$ bimodule 
$\Sigma^*\colon ={\rm Hom}_A(\Sigma,A)$, $(a\xi l)(x)=a\xi(lx)$,
for $l\in L$, $\xi \in \Sigma^*$, $a\in A$ and $x\in \Sigma$, as
\begin{equation}\label{eq:Q'}
\widetilde{Q}\colon=\{\ q\in {}_A{\rm Hom}_L(\cC,\Sigma^*)\ |\
\forall x\in \Sigma, c\in \cC \quad
c^{(1)}q(c^{(2)})(x)=q(c)(x^{^{[0]}})x^{[1]}\ \}. 
\end{equation}
\begin{proposition}\prlabel{contextmod}
Let $\Sigma$ be an $L$-$\cC$ bicomodule for a pure right coring
extension $(\cD:L)$ 
of $(\cC:A)$. Consider the corresponding Morita context
\eqref{eq:catMor}. There is a set of $k$-linear isomorphisms 
\begin{eqnarray}
&\alpha_1:&{}_L{\rm Hom}_L(\cD,T) \stackrel{\cong}{\longrightarrow} {\rm
    Nat}(V,V)^{op},\eqlabel{alpha1}\\ 
&\alpha_2:&{}^\cC{\rm End}^\cD(\cC)^{op} \stackrel{\cong}{\longrightarrow}
    {\rm Nat}(U,U)^{op},\eqlabel{alpha2}\\ 
&\alpha_3:&{}_L{\rm Hom}^\cD(\cD,\Sigma) \stackrel{\cong}{\longrightarrow}
    {\rm Nat}(V,U),\eqlabel{alpha3}\\ 
&\alpha_4:&\widetilde{Q} \stackrel{\cong}{\longrightarrow} {\rm
    Nat}(U,V),\eqlabel{alpha4} 
\end{eqnarray}
where $T$ denotes the algebra (and $L$-ring) $\mathrm{End}^\cC(\Sigma)$.
Moreover, the maps \equref{alpha1}-\equref{alpha4} establish an isomorphism of
the Morita context \eqref{eq:catMor} and the Morita context
\begin{equation}\eqlabel{eq:contextmod}
\widetilde{\mathbb M}(\Sigma)=
({}_L{\rm Hom}_L(\cD,T)\ ,\  {}^\cC{\rm End}^\cD(\cC)^{op}\ ,\
{}_L{\rm Hom}^\cD(\cD,\Sigma)\ ,\  \widetilde{Q}\ ,\newblackdiamond\
,\newdiamond ). 
\end{equation}
The algebra structures, bimodule structures and connecting
homomorphisms are given by the following formulae.
\begin{eqnarray}
&&(vv')(d)=v(d_{(1)}) v'(d_{(2)})\eqlabel{vv}\\
&&(uu')(c)=u' \big(u(c)\big)\eqlabel{uu}\\
&&(vp)(d)=v(d_{(1)})\big(p(d_{(2)})\big)\eqlabel{vp}\\
&&(pu)(d)=p(d)^{[0]}\epsilon_\cC\big(u(p(d)^{[1]})\big)\eqlabel{pu}\\
&&(qv)(c)=q(c_{[0]})v(c_{[1]})\eqlabel{qv}\\
&&(uq)(c)=q\big(u(c)\big)\eqlabel{uq}\\
&&(q\newblackdiamondd p)(c)=
c^{(1)}q({c^{(2)}}_{[0]})\big(p({c^{(2)}}_{[1]})\big)\equiv 
q(c_{[0]})\big(p({c_{[1]}})^{[0]}\big)p({c_{[1]}})^{[1]}\eqlabel{blackdiamond}
\\  
&&(p\newdiamondd q)(d)=p(d)^{[0]}q\big(p(d)^{[1]}\big)(-),\eqlabel{diamond}
\end{eqnarray}
for $v,v'\in {}_L{\rm Hom}_L(\cD,T)$, $u,u'\in {}^\cC{\rm
End}^\cD(\cC)$, $p\in {}_L{\rm Hom}^\cD(\cD,\Sigma)$, $q\in \widetilde{Q}$,
$d\in \cD$ and $c\in \cC$.
\end{proposition}
\begin{proof}
To an element $v\in {}_L{\rm Hom}_L(\cD,T)$ associate a right
$\cD$-colinear map,
$$
\Phi^v_M: {\rm Hom}^\cC(\Sigma,M)\stac{L} \cD \to {\rm
Hom}^\cC(\Sigma,M)\stac{L} \cD,\qquad \varphi_M\stac{L}d \mapsto \varphi_M\circ
v(d_{(1)})\stac{L} d_{(2)},
$$
for any right $\cC$-comodule $M$. This defines a $k$-module map
$\alpha_1:{}_L{\rm Hom}_L(\cD,T)\to {\rm Nat}(V,V)^{op}$,
$v\mapsto \Phi^v$. The bijectivity of $\alpha_1$ is proven by
constructing the inverse $\alpha_1^{-1}$, mapping $\Phi\in {\rm
Nat}(V,V)$ to the right $L$-linear map
$$
\cD\to T,\qquad d\mapsto \big((T\stac{L}\epsilon_\cD)\circ
\Phi_\Sigma\big)(1_T\stac{L} d).
$$
Since $\Sigma$ is an $L$-$\cC$ bicomodule, the map $\Sigma\to
\Sigma$, $x\mapsto lx$ is right $\cC$-colinear for any $l\in L$.
Hence $\alpha_1^{-1}(\Phi)$ is left $L$-linear by the naturality
of $\Phi$. The identity $\alpha_1^{-1}\circ \alpha_1(v)=v$, for
$v\in {}_L{\rm Hom}_L(\cD,T)$, is obvious. The other identity
$\alpha_1\circ \alpha_1^{-1}(\Phi)=\Phi$, for $\Phi\in {\rm
Nat}(V,V)$, is checked as follows. Since the right $\cD$-coaction in
$T\sstac{L} \cD$ is given by $T\sstac L\Delta_\cD$, it follows that, for all
$t\sstac L d\in T\sstac L \cD$,
\begin{equation}\label{eq:TDcoac}
(T\stac L\epsilon_\cD)\big( (t\stac L d)_{[0]}\big)\stac L (t\stac L d)_{[1]}=
t\epsilon_\cD(d_{(1)})\stac L d_{(2)}=t\stac L d.
\end{equation}
Using the right $\cD$-colinearity of $\Phi_\Sigma$ (in the second equality),
\eqref{eq:TDcoac} (in the third equality) and the right $\cC$-colinearity of
$\varphi_M$ together with the naturality of $\Phi$ (in the last equality), we
conclude that 
\begin{eqnarray*}
\big(\alpha_1\circ \alpha_1^{-1}(\Phi)\big)_M(\varphi_M\stac{L} d) &=&
\varphi_M\circ\big( (T\stac{L}\epsilon_\cD)
(\Phi_\Sigma(1_T\stac{L} d_{(1)}))\big)\stac{L} d_{(2)}\\
&=& \varphi_M\circ\big( (T\stac{L}\epsilon_\cD)\big(
\Phi_\Sigma(1_T\stac{L} d)_{[0]}\big)\big) \stac{L}
\Phi_\Sigma(1_T\stac{L} d)_{[1]}\\
&=& \big((\varphi_M\stac{L} \cD)\circ \Phi_\Sigma\big)(1_T\stac{L} d)
= \Phi_M(\varphi_M\stac{L} d),
\end{eqnarray*}
for $\varphi_M\sstac L d\in \mathrm{Hom}^\cC(\Sigma,M)\sstac L \cD$.

To an element $u\in {}^\cC{\rm End}^\cD(\cC)$ associate a map
$$
\Xi^u_M:M\to M, \qquad m\mapsto m^{[0]}(\epsilon_\cC\circ
u)(m^{[1]}),
$$
for any right $\cC$-comodule $M$. It is checked to be right
$\cD$-colinear using the relation between the $\cC$ and
$\cD$-comodule structures of $M$, the right $A$-linearity of a
$\cC$-coaction, the left $\cC$-colinearity and the right
$\cD$-colinearity of $u$ and the $\cD$-colinearity of the
$\cC$-coaction:
\begin{eqnarray*}
\big(m^{[0]}&&\!\!\!\!\!\!\!\!\!\!\!(\epsilon_\cC\circ
u)(m^{[1]})\big)_{[0]}\stac{L}
\big(m^{[0]}(\epsilon_\cC\circ u)(m^{[1]})\big)_{[1]}\\
&=& \big(m^{[0]}(\epsilon_\cC\circ u)(m^{[1]})\big)^{[0]}
\epsilon_\cC \big({\big(m^{[0]}(\epsilon_\cC\circ
u)(m^{[1]})\big)^{[1]}}_{[0]}\big)
\stac{L}{\big(m^{[0]}(\epsilon_\cC\circ u)(m^{[1]})\big)^{[1]}}_{[1]}\\
&=& m^{[0]}\epsilon_\cC \big(\big(m^{[1]}(\epsilon_\cC\circ
u)(m^{[2]})\big)_{[0]}\big)\stac{L}
\big(m^{[1]}(\epsilon_\cC\circ u)(m^{[2]})\big)_{[1]}\\
&=& m^{[0]}\epsilon_\cC
\big(u(m^{[1]})_{[0]}\big)\stac{L}u(m^{[1]})_{[1]}
={m_{[0]}}^{[0]}(\epsilon_\cC \circ u)({m_{[0]}}^{[1]})\stac{L}
m_{[1]}.
\end{eqnarray*}
This implies that we have a $k$-linear map $\alpha_2:{}^\cC{\rm
End}^\cD(\cC)^{op}\to {\rm Nat}(U,U)^{op}$, $u\mapsto \Xi^u$. In
order to prove its bijectivity, we construct the inverse
$\alpha_2^{-1}$, mapping $\Xi\in{\rm Nat}(U,U)$ to $\Xi_\cC\in{\rm
End}^\cD(\cC)$. We need to prove that $\Xi_\cC$ is left
$\cC$-colinear. Note first that, for any right $A$-module $N$ and
$n\in N$, the map $\cC\to N\otimes_A \cC$, $c\mapsto n\otimes_A c$
is right $\cC$-colinear (where $N\otimes_A \cC$ is a right
$\cC$-comodule via $N\otimes_A \Delta_\cC$). Hence, by naturality,
\begin{equation}\eqlabel{nat.Xi}
\Xi_{N\otimes_A \cC}= N\otimes_A \Xi_\cC.
\end{equation} 
On the other hand, the
coproduct in $\cC$ is right $\cC$-colinear (i.e. coassociative),
hence naturality implies $\Xi_{\cC\otimes_A \cC}\circ
\Delta_\cC=\Delta_\cC\circ \Xi_\cC$. Combining these two
observations, we have the left $\cC$-colinearity of $\Xi_\cC$
proven. The identity $\alpha_2^{-1}\circ \alpha_2(u)=u$, for
$u\in{}^\cC{\rm End}^\cD(\cC)$, follows easily by the
$\cC$-colinearity of $u$. The property $\alpha_2\circ
\alpha_2^{-1}(\Xi)=\Xi$, for $\Xi\in{\rm Nat}(U,U)$, follows 
by the commutativity of the following diagram, for any right
$\cC$-comodule $M$.
\[
\xymatrix{
M \ar[rr]^{\Xi_M} \ar[d]_{\varrho^M} && M \ar@{=}[r] \ar[d]^{\varrho^M} & M \\
M\stac A\Cc \ar[rr]_{\Xi_{M\otimes_A\Cc}=M\otimes_A\Xi_\Cc} && M\stac A\Cc
\ar[ur]_{M\otimes_A\epsilon_\Cc} 
}
\]
Commutativity of this diagram follows by \equref{nat.Xi}, the right
$\cC$-colinearity (i.e. coassociativity) of a right 
$\cC$-coaction and the naturality of $\Xi$.

To an element $p\in {}_L{\rm Hom}^\cD(\cD,\Sigma)$ associate the
right $\cD$-colinear map,
$$
\Theta^p_M: {\rm Hom}^\cC(\Sigma,M)\stac{L} \cD\to M,\qquad \varphi_M
\stac{L} d\mapsto \varphi_M \big(p(d)\big),
$$
for any right $\cC$-comodule $M$. It defines a $k$-map
$\alpha_3:{}_L{\rm Hom}^\cD(\cD,\Sigma) \to {\rm Nat}(V,U)$,
$p\mapsto \Theta^p$. We prove its bijectivity by constructing the
inverse $\alpha_3^{-1}$, mapping $\Theta\in {\rm Nat}(V,U)$ to the
right $\cD$-colinear map,
$$
\cD\to \Sigma,\qquad d\mapsto\Theta_\Sigma(1_T\stac{L} d).
$$
Its left $L$-linearity follows by the right $\cC$-colinearity of
the map $\Sigma\to \Sigma$, $x\mapsto lx$, for any $l\in L$, the
naturality of $\Theta$, and the fact that $L$ is a subalgebra of
$T$. The identity $\alpha_3^{-1}\circ \alpha_3(p)=p$, for $p\in
{}_L{\rm Hom}^\cD(\cD,\Sigma)$ is obvious, and $\alpha_3\circ
\alpha_3^{-1}(\Theta)=\Theta$, for $\Theta\in {\rm Nat}(V,U)$, follows 
by the naturality of $\Theta$, i.e. the identity 
$\varphi_M\big(\Theta_\Sigma(t\sstac L d)\big) =\Theta_M(\varphi_M\circ
t\sstac L d)$, 
for any right $\cC$-comodule $M$, $\varphi_M\in \mathrm{Hom}^\cC(\Sigma,M)$ and
$t\sstac L d\in T\sstac L\cD$.

To an element $q\in \widetilde{Q}$ associate the right $\cD$-colinear map,
\begin{eqnarray*}
\Omega^q_M: M\to{\rm Hom}^\cC(\Sigma,M)\stac{L} \cD,\qquad
m&\mapsto&
m^{[0]}q({m^{[1]}}_{[0]})(-)\stac{L}{m^{[1]}}_{[1]}\\&\equiv&
{m_{[0]}}^{[0]} q({m_{[0]}}^{[1]})(-)\stac{L}m_{[1]},
\end{eqnarray*}
for any right $\cC$-comodule $M$. The two forms of $\Omega^q_M$
are equal by the right $\cD$-colinearity of the $\cC$-coaction. It
is a well defined map by the $A$-$L$ bilinearity of $q$ and \leref{Qprop}
(2). The 
association $q\mapsto \Omega^q$ defines a $k$-map
$\alpha_4:\widetilde{Q}\to{\rm Nat}(U,V)$. We prove its bijectivity by
constructing the inverse $\alpha_4^{-1}$, mapping $\Omega\in {\rm
Nat}(U,V)$ to the right $L$ linear map
$$
\cC \to \Sigma^*, \qquad c\mapsto  \epsilon_\cC\big(({\rm
Hom}^\cC(\Sigma,\cC)\stac{L} \epsilon_\cD)(\Omega_\cC(c))(-)\big).
$$
By the right $\cC$-colinearity of the map $\cC\to
N\otimes_{A}\cC$, $c\mapsto n\otimes_{A} c$, for any right
$A$-module $N$ and $n\in N$, and the naturality of $\Omega$, the
map $\Omega_{N\sstac{A}\cC}$ can be written as a composite of $N
\otimes_{A}\Omega_{\cC}:N \otimes_{A} {\cC}\to N \otimes_{A}{\rm
Hom}^\cC(\Sigma,\cC) \otimes_L \cD$ and the obvious map $N
\otimes_{A}{\rm Hom}^\cC(\Sigma,\cC) \otimes_L
\cD\rightarrow{\rm Hom}^\cC(\Sigma,N\otimes_A \cC) \otimes_L
\cD$. Applying this fact to the case $N=A$,
we conclude on the left $A$-linearity of $\Omega_{\cC}$, hence of
$\alpha_4^{-1}(\Omega)$. 
Consider the following commutative diagram.
$$\hspace{-1.3cm}
\xymatrix { \cC
\ar[r]^{\Omega_\cC\quad\qquad}\ar[d]^{\Delta_\cC}&{\rm
Hom}^\cC(\Sigma,\cC)\stac{L} \cD\ar@<.5ex>[rr]^{\epsilon_\cC\circ
- \sstac{L} \cD}\ar[d]^{\Delta_\cC\circ - \sstac{L} \cD}&& {\rm
Hom}_A(\Sigma,A)\stac{L} \cD\ar@<.5ex>[ll]^{(-\sstac{A}\cC)\circ
\varrho^\Sigma\sstac{L} \cD} \ar[rrr]^{{\rm
Hom}_A(\Sigma,A)\sstac{L}\epsilon_\cD}
\ar[d]^{(-\sstac{A}\cC)\circ \varrho^\Sigma\sstac{L} \cD} &&& {\rm
Hom}_A(\Sigma,A)\ar[d]^{(-\sstac{A}\cC)\circ \varrho^\Sigma}
\\
\cC\stac{A}\cC\ar[r]^{\Omega_{\cC\sstac{A}\cC}\qquad\quad}
\ar[rd]_{\cC\sstac{A}\Omega_\cC}&{\rm
Hom}^\cC(\Sigma,\cC\stac{A}\cC)\stac{L}
\cD\ar@<.5ex>[rr]^{(\cC\sstac{A}\epsilon_\cC)\circ - \sstac{L}
\cD}&& {\rm Hom}_A(\Sigma,\cC)\stac{L}
\cD\ar@<.5ex>[ll]^{(-\sstac{A}\cC)\circ \varrho^\Sigma\sstac{L}
\cD} \ar[rrr]^{{\rm Hom}_A(\Sigma,\cC)\sstac{L}
\epsilon_\cD}&&&{\rm
Hom}_A(\Sigma,\cC) \\
&\cC\stac{A}{\rm Hom}^\cC(\Sigma,\cC)\stac{L} \cD
\ar@<.5ex>[rr]^{\cC\sstac{A}\epsilon_\cC\circ - \sstac{L} \cD}
\ar[u]&& \cC\stac{A}{\rm Hom}_A(\Sigma,A)\stac{L} \cD
\ar@<.5ex>[ll]^{\quad\ \cC\sstac{A}(-\sstac{A}\cC)\circ
\varrho^\Sigma\sstac{L} \cD\quad\ \ }\ar[rrr]^{\ \
\cC\sstac{A}{\rm
Hom}_A(\Sigma,A)\sstac{L}\epsilon_\cD}\ar[u]&&&\cC\stac{A}
{\rm Hom}_A(\Sigma,A)\ar[u] }
$$
The upper left square is commutative by the right
$\cC$-colinearity (i.e. coassociativity) of the coproduct in $\cC$
and the naturality of $\Omega$. The lower left triangle is
commutative by the previous observation that
$\Omega_{N\sstac{A}\cC}$ factorizes through $N\otimes_A\Omega_\cC$,
for any right $A$-module $N$. The squares in the middle column are
commutative by the isomorphism of $k$-modules ${\rm Hom}^\cC
(\Sigma,N\sstac{A}\cC)\simeq{\rm Hom}_A(\Sigma,N)$, for any right $A$-module
$N$, cf. \cite[18.10]{BrzWis:cor}. 
The upper line in the diagram gives an equivalent expression for
$\alpha_4^{-1}(\Omega)$. Comparing the incoming arrows from above and below
in $\Hom_A(\Sigma,\Cc)$ on the outer right, we conclude that
$\alpha_4^{-1}(\Omega)$ is an element of $\widetilde{Q}$. 

The identity $\alpha_4^{-1}\circ \alpha_4(q)=q$, for $q\in \widetilde{Q}$,
follows by the right $A$-linearity of $\epsilon_\cC$ and the left
$A$-linearity of $q$.  
In order to prove the converse property,
$\alpha_4\circ\alpha_4^{-1}(\Omega)=\Omega$, 
for $\Omega\in {\rm Nat}(U,V)$, consider the following diagram in $\Mm^\Dd$,
for any $M\in\Mm^\Cc$. 
$$
\xymatrix {
M \ar[rrr]^{\Omega_M\qquad\quad} \ar[rd]^{\varrho^M} \ar[ddd]^{\tau_M}
&&& {\rm Hom}^\cC(\Sigma,M)\stac{L}\cD\ar[d]^{\varrho^M\circ - \sstac{L}\cD}
\\
& M\stac{A}\cC \ar[rr]^{\Omega_{M\sstac{A}\cC}\qquad\quad}
\ar[rrd]_{M\sstac{A}\Omega_\cC} \ar[dd]_{M\ot_A\tau_\Cc}
&&{\rm Hom}^\cC(\Sigma,M\stac{A}\cC)\stac{L}\cD 
\\
&&& M\stac{A}{\rm Hom}^\cC(\Sigma,\cC)\stac{L}\cD\ar[u]
\\
M\stac L\Dd \ar[r]^-{\rho^M\ot_L\Dd}
& M\stac A\Cc\stac L\Dd \ar[rr]^-{M\ot_A\Omega_\Cc\ot_L\Dd}
&&M\stac A\Hom^\Cc(\Sigma,\Cc)\stac L\Dd\stac L\Dd
\ar@<-.5ex>[u]_{M\ot_A\Hom^\Cc(\Sigma,\Cc)\ot_L\epsilon_\Dd\ot_L\Dd} 
}
$$
The commutativity of the upper quadrangle follows by the naturality of
$\Omega_M$. The left lower quadrangle commutes by the right $\Dd$-colinearity
of $\rho^M$ and the right lower quadrangle does by the right $\Dd$-colinearity
of $\Omega_\Cc$ and the fact that $\epsilon_\cD$ is the counit of $\Dd$.
The property that $\Omega_{M\ot_A\Cc}$ factorizes trough $M\ot\Omega_\Cc$ implies
that the triangle commutes as well. Evaluating the upper and lower
incoming arrows in $\Hom^\Cc(\Sigma,M\ot_A\Cc)\ot_L\Dd$, on an element $m\in
M$, one finds 
\begin{eqnarray*}
(\rho^M\circ -)\circ\Omega_M(m)&=&
  {m_{[0]}}^{[0]}\ot_A(\Hom^\Cc(\Sigma,\Cc)\ot_A\epsilon_\Dd)
(\Omega_\Cc({m_{[0]}}^{[1]}))\ot_Lm_{[1]}.  
\end{eqnarray*}
Application of  $(M\ot_A\epsilon_\Cc)\circ -\ot_L\Dd$ to both sides of this
equation yields the required identity 
$\Omega_M=\alpha_4\circ\alpha_4^{-1}(\Omega)_M$ in
$\Hom^\Cc(\Sigma,M)\ot_L\Dd$. 

The proof is completed by showing that the maps $\alpha_1$,
$\alpha_2$, $\alpha_3$ and $\alpha_4$ define a morphism of Morita
contexts. Indeed, it is straightforward to check that, for
$v,v'\in{}_L{\rm Hom}_L(\cD,T)$, $u,u'\in {}^\cC{\rm End}^\cD(\cC)$,
$p\in {}_L{\rm Hom}^\cD(\cD,\Sigma)$ and $q\in \widetilde{Q}$,
\begin{eqnarray*}
&\alpha_1(v)\circ \alpha_1(v')=\alpha_1(v'v)\qquad
&\alpha_2(u)\circ \alpha_2(u')=\alpha_2(u'u)\\
&\alpha_3(p)\circ \alpha_1(v)=\alpha_3(v\ p)\qquad
&\alpha_2(u)\circ \alpha_3(p)=\alpha_3(p\ u)\\
&\alpha_1(v)\circ \alpha_4(q)=\alpha_4(q\ v)\qquad
&\alpha_4(q)\circ \alpha_2(u)=\alpha_4(u\ q)\\
&\alpha_3(p)\circ \alpha_4(q)=\alpha_2(q\newblackdiamondd p)\qquad
&\alpha_4(q)\circ \alpha_3(p)=\alpha_1(p\newdiamondd q).
\end{eqnarray*}
This ends the proof.
\end{proof}

\begin{remark}\label{rem:k}
The commutative base ring $k$ can be considered as a trivial ($k$-) coring,
which is a pure right extension of any $A$-coring $\cC$. A right $\cC$-comodule
$\Sigma$ can be viewed as a $k$-$\cC$ bicomodule for the pure right coring
extension $k$ of $\cC$, hence there is an associated Morita context
$\widetilde{\mathbb M}(\Sigma)=
\big(\ {\Hom_k(k,T)}, {^\Cc\End}(\Cc)^{\rm op}, {\Hom_k}(k,\Sigma),
\widetilde{Q}\equiv Q',\newblackdiamond,\newdiamond\ \big)$, as in
\equref{eq:contextmod}. Obviously, 
${\Hom_k(k,T)}\cong T$ and ${\Hom_k}(k,\Sigma)\cong \Sigma$. By a hom-tensor
relation for comodules \cite[18.11]{BrzWis:cor}, 
${^\Cc\End}(\Cc)^{\rm op}\cong \*C$. 
By \leref{Qprop} (1), $Q\cong Q'\equiv \widetilde{Q}$. Composing these
isomorphisms with the Morita maps in \equref{eq:contextmod},
one obtains 
the structure maps of the Morita context \eqref{eq:SiMor}.
Hence the Morita context
$\widetilde{\mathbb M}(\Sigma)$, associated to $\Sigma$ as a
$k$-$\cC$ bicomodule, coincides with ${\mathbb M}(\Sigma)$, associated to a 
right $\cC$-comodule $\Sigma$ in \eqref{eq:SiMor}.
\end{remark}

\begin{lemma} \label{lem:jids} \lelabel{jids}
Let the $L$-coring $\cD$ be a pure right extension of the
$A$-coring $\cC$ 
and let $\Sigma$ be an object in ${}_L\cM^\cC$. Consider the Morita context
$\widetilde{\mathbb M}(\Sigma)$ in \equref{eq:contextmod}. If there
exist finite sets of elements 
$\{j_\ell\}\subset{_L\Hom^\Dd}(\Dd,\Sigma)$ and 
$\{\tildej_\ell\}\subset \widetilde{Q}$ such
that 
$\sum_\ell\tildej_\ell\newblackdiamondd j_\ell=\Cc$ (i.e. connecting map
$\newblackdiamond$ in \equref{blackdiamond} is surjective), then 
\begin{zlist}
\item the identity
$
\sum_\ell\tildej_\ell(c_{[0]})\big( j_\ell(c_{[1]})\big)=\epsilon_\Cc(c)
$ holds, for all $c\in \cC$;
\item the identity
$\sum_\ell {m_{[0]}}^{[0]}
  \tildej_\ell({m_{[0]}}^{[1]})\big(j_\ell(m_{[1]})\big) =m$ holds, for any
  right $\cC$-comodule $M$ and $m\in M$. 
\end{zlist}
\end{lemma}
\begin{proof}
{\underline {(1)}} 
This follows by applying $\epsilon_\Cc$ to \equref{blackdiamond}.

{\underline {(2)}} 
By the $\cD$-colinearity of a right $\cC$-coaction in a right $\cC$-comodule 
$M$, for any $m\in M$,
\[
\sum_\ell {m_{[0]}}^{[0]}
\tildej_\ell({m_{[0]}}^{[1]})\big(j_\ell(m_{[1]})\big) 
=\sum_\ell {m}^{[0]}
\tildej_\ell({m^{[1]}}_{[0]})\big(j_\ell({m^{[1]}}_{[1]})\big) 
=m^{[0]} \epsilon_\cC(m^{[1]})=m,
\]
where we used part (1) in the second equation.
\end{proof}

\begin{proposition}\label{prop:SiAmod}
Let the $L$-coring $\cD$ be a pure right extension of the
$A$-coring $\cC$. 
Take $\Sigma\in {}_L\cM^\cC$ and consider the associated Morita context
$\widetilde{\mathbb M}(\Sigma)$ in \equref{eq:contextmod}. 
If both the map $\newblackdiamond$ in \equref{blackdiamond} and
the counit $\epsilon_\cC$ are surjective then $\Sigma$ is a generator of
right $A$-modules.
\end{proposition}
\begin{proof}
Choose elements $\{j_\ell\}\subset{_L\Hom^\Dd}(\Dd,\Sigma)$ and 
$\{\tildej_\ell\}\subset \widetilde{Q}$ such that 
$\sum_\ell\tildej_\ell\newblackdiamondd j_\ell=\Cc$, and an element $c\in \cC$
such that $\epsilon_\cC(c)=1_A$. Fix finite sets $\{x_i\}\subset \Sigma$ and
$\{\xi_i\}\subset \Sigma^*$ 
such that
$$
 \sum_i \xi_i\stac L x_i=\sum_\ell\tildej_\ell(c_{[0]})\stac{L}
j_\ell(c_{[1]}). 
$$ 
Then $\sum_i \xi_i(x_i)=\epsilon_\Cc(c)=1_A$, 
by \leref{jids} (1), which proves the claim.
\end{proof}

A lemma by Beck (cf. a dual version of 
\cite[3.3 Proposition 3]{BarrWells:ttt}) states that  
if $F:{\mathcal A}\to {\mathcal B}$  is a comonadic functor and $\lambda$ is a
split epimorphism in the category ${\mathcal A}$, such that $F(\lambda)$ has a
kernel in ${\mathcal B}$, then also $\lambda$ has a kernel in ${\mathcal A}$
and $F$ preserves this kernel.
Since for an $A$-coring $\cC$ the forgetful functor $\cM^\cC\to \cM_A$ is
obviously comonadic, the next lemma follows by this general result. 
Still, for the convenience of the reader we include a complete proof, which is
a simplified version of
Beck's arguments (as in our case the target category $\cM_A$ is Abelian). 
\begin{lemma}\lelabel{summand}
The following statements about right comodules $M$ and $N$ for an $A$-coring
$\cC$ are equivalent.
\begin{enumerate}[(i)]
\item There exists an index set $\{\, \ell\, \}$ of finite order $s$ and
  collections of morphisms $\{\kappa_\ell\}\subset \mathrm{Hom}^\cC(M,N)$ and 
$\{\lambda_\ell\}\subset \mathrm{Hom}^\cC(N,M)$, indexed by $\ell=1\dots s$,
  such that $\sum_\ell \lambda_\ell \circ\kappa_\ell=M$;  
\item $M$ is a direct summand of the direct sum $N^{s}$ as a right
  $\Cc$-comodule. 
\end{enumerate}
\end{lemma}

\begin{proof}
In terms of $\{ \kappa_\ell\}$ and $\{\lambda_\ell\}$ as in (i),
construct a  map  
$\kappa:M\to N^s$ defined by $\kappa(m)=(\kappa_\ell(m))_\ell$. 
Clearly, $\kappa$ has a left inverse, $\lambda:N^s \to M$,
$\lambda((n_\ell)_\ell) =\sum_\ell \lambda_\ell(n_\ell)$. Recall
that the category of comodules of any $A$-coring $\cC$ has direct sums that
do coincide with the direct sums in $\Mm_A$ (see e.g. \cite[3.3]{Wis:catcom}), 
so $\lambda$ and $\kappa$ are morphisms in $\Mm^\Cc$.
In particular, $\lambda$ and $\kappa$ are right $A$-linear, so 
there is a split exact sequence in the Abelian category of right $A$-modules,
$$
\xymatrix{
0\ar[r]&M'\ar[r]^{\nu}&N^s\ar[r]^{\lambda}&M\ar[r]&0\ ,
}
$$
where $(M',\nu)$ is the kernel of $\lambda$ in $\cM_A$.
The proof is completed by showing that $M'$ is a right $\Cc$-comodule and
$\nu$ is a $\cC$-colinear section. 
Denoting a right $A$-linear retraction of $\nu$ by $\varpi$, introduce a right
$A$-module map
$$
\varrho^{M'}\colon =(\varpi\stac A \cC)\circ \varrho^{N^s}\circ \nu:\quad
M'\to M'\stac A \cC.
$$
Since the $A$-module maps satisfy $\nu\circ \varpi +\kappa\circ\lambda=N^s$,
the colinearity of $\kappa$ 
and $\lambda$ implies that $\nu\circ \varpi$ is $\cC$-colinear. Hence
\begin{equation}\label{eq:kcolin}
(\nu \stac A \cC)\circ \varrho^{M'}=(\nu \circ \varpi\stac A\cC)\circ 
\varrho^{N^s}\circ \nu= \varrho^{N^s}\circ \nu.
\end{equation}
Composing \eqref{eq:kcolin} by $\varpi$ on both sides and using the
colinearity of $\nu\circ \varpi$, one obtains
also $\varrho^{M'}\circ \varpi=(\varpi\sstac A \cC)\circ \varrho^{N^s}$.
Furthermore, using \eqref{eq:kcolin} (in the first, second and fourth
equalities) and the coassociativity of $ \varrho^{N^s}$ (in the third
equality), one checks that 
\begin{eqnarray*}
(\nu \stac A\cC\stac A\cC)\circ (\varrho^{M'}\stac A \cC)\circ \varrho^{M'}
&=& (\varrho^{N^s}\circ \nu \stac A \cC)\circ \varrho^{M'}
= (\varrho^{N^s}\stac A \cC)\circ \varrho^{N^s}\circ \nu \\
&=& (N^s\stac A\Delta_\cC)\circ \varrho^{N^s}\circ \nu
=(N^s\stac A\Delta_\cC)\circ (\nu \stac A \cC)\circ \varrho^{M'}\\
&=&(\nu \stac A\cC\stac A\cC)\circ (M'\stac A\Delta_\cC)\circ\varrho^{M'}.
\end{eqnarray*}
Since $\nu \sstac A\cC\sstac A\cC$ is a (split) monomorphism, this implies the
coassociativity of $\varrho^{M'}$. Finally, by the counitality of 
$\varrho^{N^s}$, 
$$
(M'\stac A \epsilon_\cC)\circ \varrho^{M'}=(M'\stac A \epsilon_\cC)\circ
(\varpi\stac A \cC)\circ \varrho^{N^s}\circ \nu  =\varpi \circ \nu=M',
$$
which finishes the proof of the implication $(i)\Rightarrow (ii)$. 
The converse implication is obvious.
\end{proof}

Since $L$-$\cC$ bicomodules, for an algebra $L$ and an $A$-coring $\cC$, can
be identified with right comodules for the $L^{op}\sstac k A$-coring
$L^{op}\sstac k \cC$, \leref{summand} applies also to $L$-$\cC$ bicomodules. 

\begin{theorem}\thlabel{surjectivity}
Let the $L$-coring $\cD$ be a pure right extension of the
$A$-coring $\cC$.  
Take $\Sigma\in {}_L\cM^\cC$ and consider the Morita context
$\widetilde{\mathbb M}(\Sigma)$ in \equref{eq:contextmod}. In particular, put
$T\colon = \mathrm{End}^\cC(\Sigma)$.
\begin{enumerate}
\item 
The map $\newblackdiamond$ in \equref{blackdiamond} is surjective if and only
if $\Sigma$ is a Galois 
$\cC$-comodule and $\Sigma$ is a direct summand of a direct sum
$(T\otimes_L\Dd)^s$ as 
$T$-$\Dd$-bicomodule, for an appropriate finite integer $s$.
\item 
The Morita context $\widetilde{\mathbb M}(\Sigma)$ is strict if and only if
$\Sigma$ is a Galois  $\cC$-comodule,
$\Sigma$ is a direct summand of $(T\otimes_L\Dd)^{s}$ and $T\otimes_L\Dd$ 
is a direct summand of $\Sigma^{z}$, both as $T$-$\Dd$-bicomodules, where $s$ 
and $z$ are finite integers. 
\end{enumerate}
\end{theorem}

\begin{proof}
\ul{(1)}
By surjectivity of $\newblackdiamond$, there exist elements
$\{j_\ell\}\subset {_L\Hom^\Dd}(\Dd,\Sigma)$ and 
$\{\tildej_\ell\}\subset \widetilde{Q}$ (cf. \eqref{eq:Q'}) such that 
$\sum_\ell\tildej_\ell\newblackdiamondd j_\ell=\Cc$. In terms of the
isomorphisms $\alpha_3$ and $\alpha_4$ in \prref{contextmod}, we denote
$\alpha_3(j_\ell)=J_\ell$ and $\alpha_4(\tildej_\ell)=\tildeJ_\ell$, for all
values of $\ell$.  

First we check that the surjectivity of $\newblackdiamond$ implies the Galois
property of the right $\cC$-comodule $\Sigma$. To this end, we
construct the inverse of the canonical natural transformation
\eqref{eq:can}. For any right $A$-module $N$, put
$$
\Upsilon_N : N\stac{A} \cC\to \Hom_A(\Sigma,N)\stac{T} \Sigma, \qquad
n\stac{A} c \mapsto \sum_\ell n\tildej_\ell(c_{[0]})(-)\stac{T}
j_\ell(c_{[1]}). 
$$
Since $T$ is an $L$-ring, $\Upsilon_N$ is a well defined map.
It is natural in $N$, being a sum of composites of natural morphisms,
\[
\xymatrix{
N\stac{A}\Cc \ar[rr]^-{(\tildeJ_\ell)_{N\otimes_{A}\Cc}} &&
\stackrel{\displaystyle \quad \Hom^\Cc\,(\Sigma,N\stac{A}\Cc)\, \stac{L}
  \Dd} {\cong\Hom_A(\Sigma,N)\stac{T}T\stac{L} \Dd} 
\ar[rrr]^-{\Hom_A(\Sigma,N)\otimes_{T} (J_\ell)_\Sigma} &&&
\Hom_A(\Sigma,N)\stac{T} \Sigma\ .
}
\]
We claim that $\Upsilon$ 
yields the inverse of the canonical natural transformation
\eqref{eq:can}. Indeed, we find that, for $n\sstac{A} c\in N\sstac{A} \cC$,
\begin{eqnarray*}
\can_N(\Upsilon_N(n \stac{A} c))&=&
\sum_\ell n\stac{A} {\tildej_\ell}(c_{[0]})\big(j_\ell(c_{[1]})^{[0]}\big)
j_\ell(c_{[1]})^{[1]}\\ 
&=&n\stac{A} \sum_\ell (\tildej_\ell\newblackdiamondd j_\ell)(c) =
n\stac{A} c.
\end{eqnarray*}
On the other hand, for $\phi_N\sstac{T} x\in \Hom_A(\Sigma,N)\sstac{T} \Sigma$,
\begin{eqnarray*}
\Upsilon_N\big(\can_N(\phi_N\stac{T} x)\big)&=&
\sum_\ell\phi_N(x^{[0]}) {\tildej_\ell}({x^{[1]}}_{[0]})(-)\stac{T}
j_\ell({x^{[1]}}_{[1]})\\ 
&=&\sum_\ell \phi_N\big( {x_{[0]}}^{[0]}{\tildej_\ell}({x_{[0]}}^{[1]})(-)\big)
\stac{T} j_\ell({x_{[1]}}),
\end{eqnarray*} 
where in the second  equality we used the right $\cD$-colinearity of 
the $\cC$-coaction in $\Sigma$ and the right $A$-linearity of
$\phi_N$. 
Using \leref{3.4} (2) at the right $\cC$-comodule $\Sigma$, we
conclude that $x^{[0]} {\tildej_\ell}({x^{[1]}})(-)$ is an element of $T$, for
all $x\in \Sigma$ and any value of the index $\ell$. Hence
$$
\sum_\ell \phi_N\big(
{x_{[0]}}^{[0]}{\tildej_\ell}({x_{[0]}}^{[1]})(-)\big)\stac{T}  
j_\ell({x_{[1]}})=\phi_N\stac{T} 
\sum_\ell
{x_{[0]}}^{[0]}{\tildej_\ell}({x_{[0]}}^{[1]})\big(j_\ell({x_{[1]}})\big) = 
\phi_N\stac{T} x,
$$ 
where the last equality follows by Lemma \ref{lem:jids} (2), applied to the
right $\cC$-comodule $\Sigma$. This
proves that if the Morita map $\newblackdiamond$ is surjective then 
$\Sigma$ is a Galois $\cC$-comodule. 

Next we prove that $\Sigma$ is a direct summand of $(T\otimes_L\Dd)^s$ as
$T$-$\Dd$-bicomodule, where $s$ is the cardinality of the index set $\{\ell\}$
above. For any value of $\ell$, put
\begin{eqnarray}\label{eq:kappa}
&\kappa_\ell\colon =(\widetilde{J}_\ell)_\Sigma\ :\ \Sigma\to T\stac{L}
  \cD,\qquad  
&x\mapsto {x_{[0]}}^{[0]}{\tildej_\ell}({x_{[0]}}^{[1]})(-)\stac{L}
x_{[1]},\quad \textrm{and}\\
&{\widetilde \kappa}_\ell\colon = ({J}_\ell)_\Sigma\ :\ T\stac{L} \cD\to
  \Sigma,\qquad 
&t\stac{L} d \mapsto t\big(j_\ell(d)\big).\nonumber
\end{eqnarray}
Since the left $L$-action in $\Sigma$ is right $\cC$-colinear by assumption,
and any $\cC$-colinear map is $\cD$-colinear, both the right $\cC$, and
$\cD$-coactions in $\Sigma$ are left $L$-linear. This way $\kappa_\ell$, being
a composite of left $L$-linear right $\cD$-colinear maps, is left $L$-linear
and right $\cD$-colinear. The map ${\widetilde \kappa}_\ell$ is obviously left
$L$-linear and it
is right $\cD$-colinear by the colinearity of $j_\ell$ and any $t\in
T$. It follows by 
Lemma \ref{lem:jids} (2) that $\sum_\ell{\widetilde
  \kappa}_\ell\circ\kappa_\ell(x)=x$. By \leref{summand}, this implies 
that $\Sigma$ is a direct summand of $(T\otimes_L\Dd)^s$ as
$T$-$\Dd$-bicomodule, where $s$ is the cardinality of the index set
$\{\ell\}$.  

In order to prove the converse statement, we make use of \leref{summand}
again. Since $\Sigma$ is a direct summand of $(T\otimes_L\Dd)^s$ as 
$T$-$\Dd$-bicomodule, there exist maps
$\kappa_\ell\in{_T\Hom^\Dd}(\Sigma,T\sstac{L}\Dd)$ and ${\widetilde
  \kappa}_\ell\in{_T\Hom^\Dd}(T\sstac{L}\Dd,\Sigma)$, for $\ell=1\dots s$, such
that 
$\sum_\ell{\widetilde \kappa}_\ell\circ\kappa_\ell=\Sigma$. We can define
maps $j_\ell$ and $\tildej_\ell$ as follows. For any value of $\ell$, put
\begin{equation}\eqlabel{defj}
j_\ell: \cD\to \Sigma, \qquad 
d\mapsto {\widetilde \kappa}_\ell(1_T\stac{L} d).
\end{equation}
By the colinearity of ${\widetilde \kappa}_\ell$, $j_\ell$ is right
$\cD$-colinear. 
Since ${\widetilde \kappa}_\ell$ is left $T$-linear and $L$ is a subalgebra of
$T$, $j_\ell$ is 
also left $L$-linear.
Furthermore, $\Sigma$ is a Galois $\cC$-comodule by assumption, hence we
can set
\begin{equation}\eqlabel{deftildej}
{\tildej}_\ell\colon =[\Sigma^*\stac{T} 
(T\stac{L} \epsilon_\cD)\circ \kappa_\ell]\circ \can_A^{-1}: \quad
\cC \to \Sigma^*,
\end{equation}
for $\ell =1\dots s$. The map $\can_A$ is left $A$-linear and right
$\cC$-colinear. Hence it is also 
right $\cD$-colinear, so, in particular, right $L$-linear.
Since this way ${\tildej}_\ell$ is a composite of left $A$-linear, right
$L$-linear maps, it is left $A$-linear and right $L$-linear. Let us prove that
${\tildej}_\ell\in \widetilde{Q}$. Consider the following commutative diagram. 
\[
\xymatrix{
\Cc \ar[r]^-{\can_A^{-1}} \ar[d]_{\Delta_\Cc} 
& \Hom_A(\Sigma,A)\stac T\Sigma
\ar[rrrr]^-{\Hom_A(\Sigma,A)\otimes_T(T\otimes_L\epsilon_\Dd)\circ
\kappa_\ell} 
\ar[d]^{(-\otimes_A\Cc)\circ\varrho^\Sigma\otimes_T\Sigma} 
&&&& \Hom_A(\Sigma,A) \ar[d]^{(-\otimes_A\Cc)\circ\varrho^\Sigma} \\
\Cc\stac A\Cc \ar[r]^-{\can_\Cc^{-1}}
& \Hom_A(\Sigma,\Cc)\stac T\Sigma
\ar[rrrr]^-{\Hom_A(\Sigma,\Cc)\otimes_T(T\otimes_L\epsilon_\Dd)\circ
\kappa_\ell} 
&&&& \Hom_A(\Sigma,\Cc) \\
\Cc \ar[r]_-{\can_A^{-1}} \ar[u]^{f_c\otimes_A\Cc} 
& \Hom_A(\Sigma,A)\stac T\Sigma \ar[rrrr]_-{\Hom_A(\Sigma,A)\otimes_T
(T\otimes_L\epsilon_\Dd)\circ \kappa_\ell} 
\ar[u]_{(f_c\circ-)\otimes_T\Sigma}
&&&& \Hom_A(\Sigma,A) \ar[u]_{f_c\circ-}
}
\]
Here $f_c$ denotes a right $A$-linear morphism from $A$ to $\Cc$, 
defined for an element $c\in\Cc$ as $f_c(a)=ca$. The lower square on the left
hand side 
commutes because of the naturality of $\can^{-1}$. By the explicit form
\eqref{eq:can} of $\mathrm{can}$, for $\xi\sstac T x\in \Sigma^*\sstac T
\Sigma$, 
$$
\mathrm{can}_\cC\circ [(-\stac A \cC)\circ \rho^\Sigma \stac T \Sigma]
(\xi \stac T x)
=\xi(x^{[0]})x^{[1]}\stac A x^{[2]}
=\Delta_\cC\circ \mathrm{can}_A(\xi \stac T x).
$$
Hence also the upper square on the left hand side commutes.
The commutativity of the squares on the right hand side is obvious.
Both the upper and lower horizontal lines express the map ${\tildej}_\ell$.  
In terms of the map $f_c$ introduced above, $\Delta_\cC(c)=(f_{c^{(1)}}\sstac A
\cC)(c^{(2)})$. Hence the commutativity of the diagram
implies, for all $c\in \cC$ and $x\in \Sigma$,  
\begin{eqnarray*}
{\tildej}_\ell(c)(x^{[0]})x^{[1]}
&=&\big([\mathrm{Hom}_A(\Sigma,\cC)\stac T (T\stac L \epsilon_\cD)\circ
  \kappa_\ell] \circ \mathrm{can}_\cC^{-1}\circ \Delta_\cC(c)\big)(x)\\
&=&\big([\mathrm{Hom}_A(\Sigma,\cC)\stac T (T\stac L \epsilon_\cD)\circ
  \kappa_\ell] \circ \mathrm{can}_\cC^{-1}\circ(f_{c^{(1)}}\stac A
\cC)(c^{(2)})\big)(x)\\
&=&c^{(1)}{\tildej}_\ell(c^{(2)})(x).
\end{eqnarray*}
That is, $\tildej_\ell\in \widetilde{Q}$, for all values of the index $\ell$.
The surjectivity of (the ${}^\cC\mathrm{End}^\cD(\cC)$-${}^\cC
\mathrm{End}^\cD(\cC)$ bilinear map) $\newblackdiamond$ is proven by showing
$\sum_\ell{\tildej}_\ell\newblackdiamondd j_\ell=\Cc$.
Use the right
$\cD$-colinearity of $\can_A^{-1}$ (in the second equality), the right 
$\cD$-colinearity of $\kappa_\ell$ (in the third one) and the left
$T$-linearity of $\widetilde{\kappa}_\ell$ (in the penultimate one) to 
compute the composite of the right $\cD$-coaction $\tau_\cC$ on $\cC$ with 
$\sum_\ell {\tildej}_\ell  \sstac{L} j_\ell$. It yields
\begin{eqnarray}
\sum_\ell ({\tildej}_\ell  \stac{L} j_\ell)\circ \tau_\cC
&=&\sum_\ell\{[\Sigma^*\stac{T} (T\stac{L} \epsilon_\cD)\circ
  \kappa_\ell] \stac{L} j_\ell\}\circ (\can_A^{-1}\stac{L} \cD)\circ \tau_\cC
\nonumber\\
&=& \sum_\ell\{\Sigma^*\stac{T} 
[(T\stac{L} \epsilon_\cD)\circ \kappa_\ell \stac{L} j_\ell]\circ \tau_\Sigma\}
\circ \can_A^{-1}\nonumber\\
&=& \sum_\ell[\Sigma^*\stac{T} 
(T\stac{L} \epsilon_\cD\stac{L} j_\ell) \circ 
(T\stac{L} \Delta_\cD)\circ \kappa_\ell]\circ \can_A^{-1}\nonumber\\
&=& \sum_\ell [\Sigma^*\stac{T} (T\stac{L} j_\ell)\circ \kappa_\ell]
\circ \can_A^{-1}.\nonumber\\
&=& [\Sigma^*\stac{T}\sum_\ell {\widetilde \kappa}_\ell\circ
  \kappa_\ell] 
\circ \can_A^{-1} = \can_A^{-1}.
\label{eq:jtj}
\end{eqnarray}
Note that the evaluation map $\Sigma^*\otimes_T\Sigma\to A$,
$\xi\sstac T x\mapsto \xi(x)$ 
is equal to $\epsilon_\cC\circ \can_A$. Hence equation \eqref{eq:jtj} implies 
\begin{eqnarray*}
\sum_\ell\tildej_\ell\newblackdiamondd j_\ell
&=&\sum_\ell(\epsilon_\cC\circ
\can_A\stac{A}\cC)\circ(\Sigma^*\stac{T} \rho^\Sigma)\circ
({\tildej_\ell}\stac{L} j_\ell)\circ \tau_\cC\\ 
&=&(\epsilon_\cC\circ \can_A \stac{A}\cC)\circ(\Sigma^*\stac{T}
\rho^\Sigma)\circ \can_A^{-1}
=(\epsilon_\cC \stac{A} \cC)\circ \Delta_\cC
=\cC,
\end{eqnarray*}
where the third equality follows by the right $\cC$-colinearity
of $\can_A$. 

\ul{(2)} 
Suppose that the Morita context \equref{eq:contextmod} is
strict. In view of part (1), we have to prove only that $T\sstac{L}\Dd$ is a
direct summand of $\Sigma^z$, for some integer $z$. By the surjectivity of
$\newdiamond$, there exist elements $\{h_i\}\subset{_L\Hom^\Dd}(\Dd,\Sigma)$
and $\{{\widetilde h}_i\}\subset \widetilde{Q}$ such that  
$\sum_i (h_i\newdiamondd {\widetilde h}_i)(d)=\epsilon_\cD(d) 1_T$, for $d\in
\cD$. Similarly to \eqref{eq:kappa}, for any value of $i$, we define 
left $L$-linear and right $\Dd$-colinear morphisms,  
\begin{eqnarray*}
&\lambda_i\ :\ \Sigma\to T\stac{L}
  \cD,\qquad  
&x\mapsto {x_{[0]}}^{[0]}{\widetilde h}_i({x_{[0]}}^{[1]})(-)\stac{L}
x_{[1]}\quad \textrm{and}\\
&{\widetilde \lambda}_i\ :\ T\stac{L} \cD\to
  \Sigma,\qquad 
&t\stac{L} d \mapsto t\big(h_i(d)\big).
\end{eqnarray*}
They satisfy, for any $t\sstac{L} d\in T\sstac{L} \cD$,
\begin{eqnarray*}
\sum_i \lambda_i \big({\widetilde \lambda}_i(t\stac{L} d)\big)&=&
\sum_i{t(h_i(d))_{[0]}}^{[0]}{\widetilde
  h}_i\big({t(h_i(d))_{[0]}}^{[1]}\big)(-) 
\stac{L} {t(h_i(d))_{[1]}}\\
&=&\sum_it\big(h_i(d_{(1)})^{[0]}{\widetilde
  h}_i(h_i(d_{(1)})^{[1]})(-)\big) \stac{L} d_{(2)}\\
&=&\sum_i t\big((h_i\newdiamondd \widetilde{h}_i)(d_{(1)})\big)\stac L
d_{(2)} 
=t\epsilon_\cD(d_{(1)})\stac{L} d_{(2)}
=t\stac{L} d,
\end{eqnarray*}
where the second 
 equality follows by $\cC$- and $\cD$-colinearity of $t\in T$,
 $\cD$-colinearity of $h_i$,
and right $A$-linearity of $t\in T$. 
By \leref{summand}, we conclude that $T\sstac{L}\Dd$ is 
a direct summand of $\Sigma^z$, where $z$ is the cardinality of the index set
$\{i\}$. 

Finally, we show that if $\Sigma$ is a Galois $\cC$-comodule and
$T\sstac{L}\Dd$ is a direct summand of $\Sigma^z$, 
then $\newdiamond$ is surjective. Let
$\{{\widetilde \lambda}_i\}\subset {}_T\mathrm{Hom}^\cD(T\sstac L \cD,\Sigma)$
and $\{\lambda_i\}\subset {}_T\mathrm{Hom}^\cD(\Sigma,T\sstac L \cD)$ be sets
of morphisms, such 
that $\sum_i\lambda_i\circ{\widetilde \lambda}_i=T\sstac{L}\Dd$. Repeating the
arguments in part (1) (cf. \equref{defj} and \equref{deftildej}), we define
maps $h_i\in{_L\Hom^\Dd}(\Dd,\Sigma)$ and ${\widetilde h}_i\in \widetilde{Q}$
as
$$  
h_i\colon ={\widetilde \lambda}_i(1_T\stac{L} -)\qquad \textrm{and}\qquad 
{\widetilde h}_i\colon=[\Sigma^*\stac{T} 
(T\stac{L} \epsilon_\cD)\circ \lambda_i]\circ \can_A^{-1}.
$$
For any $x\in \Sigma$, the association $a\mapsto xa$ defines a right $A$-module
map $A\to \Sigma$.
By naturality of the canonical maps \eqref{eq:can}, for $x\in \Sigma$,
$c\in \cC$, and any value of the index $i$,
\begin{eqnarray}\label{eq:xht}
x{\widetilde h}_i(c)(-)
&=&x\{[{\rm Hom}_A(\Sigma,A)\stac{T} (T\stac{L} \epsilon_\cD)\circ \lambda_i]
\circ {\rm can}_A^{-1}(c)\}(-)   \\
&=&[{\rm End}_A(\Sigma)\stac{T} (T\stac{L} \epsilon_\cD)\circ \lambda_i]
\circ {\rm can}_\Sigma^{-1}(x\stac{A} c).\nonumber
\end{eqnarray}
By \eqref{eq:can}, $\can_\Sigma(\Sigma\sstac T x)=x^{[0]}\sstac A
x^{[1]}$, for $x\in \Sigma$. Hence \eqref{eq:xht} implies the following
equality of right $A$-linear endomorphisms of $\Sigma$.
$$
x^{[0]}{\widetilde h}_i(x^{[1]})(-)=(T\stac{L} \epsilon_\cD)\circ
\lambda_i(x),
$$
for all $x\in \Sigma$ and any value of the index $i$.
Then it follows that, for $d\in \cD$,
\begin{eqnarray*}
\sum_i (h_i\newdiamondd {\widetilde h}_i)(d)
&=&
\sum_i h_i(d)^{[0]}{\widetilde h}_i\big(h_i(d)^{[1]}\big)(-)
=\sum_i(T\stac{L} \epsilon_\cD)\big(\lambda_i(h_i(d))\big)\\
&=&(T\stac{L} \epsilon _\cD)\big(\sum_i (\lambda_i\circ {\widetilde 
\lambda}_i)(1_T\stac{L} d)\big) 
=\epsilon_\cD(d)1_T.
\end{eqnarray*}
This proves that $\sum_i h_i\newdiamondd {\widetilde h}_i$ is the
unit element of the algebra ${}_L\mathrm{Hom}_L(\cD,T)$, and thus the
surjectivity of the (${}_L\mathrm{Hom}_L(\cD,T)$-${}_L\mathrm{Hom}_L(\cD,T)$
bilinear) map $\newdiamond$. 
\end{proof}

\begin{remark}\relabel{numbers}
It follows by the proof of \thref{surjectivity} that the finite number $s$
in the theorem can be chosen equal to the cardinality of the sets
$\{j_\ell\}\subset{_L\Hom^\Dd}(\cD,\Sigma)$ and 
$\{\tildej_\ell\}\subset \widetilde{Q}$, such that 
$\sum_\ell \tildej_\ell\newblackdiamondd j_\ell=\cC$. Similarly, the number
$z$ can be chosen equal to the cardinality of the sets 
$\{h_i \}\subset{_L\Hom^\Dd}(\cD,\Sigma)$ and 
$\{\widetilde{h}_i\}\subset \widetilde{Q}$, such that 
$\sum_i  h_i \newdiamondd \widetilde{h}_i=\epsilon_\cD(-)1_T$.  
\end{remark}

\begin{proposition}\label{prop:eqproj}\prlabel{eqproj}
Let the $L$-coring $\Dd$ be a pure right extension of the
$A$-coring $\Cc$ and 
$\Sigma\in{_L\Mm^\Cc}$. Consider the associated Morita context
$\widetilde{\mathbb M}(\Sigma)$ in
\equref{eq:contextmod}. 
If the connecting map $\newblackdiamond$ in \equref{blackdiamond} is
surjective, 
then $\Sigma$ is a $\cD$-equivariantly $L$-relative projective left
module of the algebra $T=\mathrm{End}^\cC(\Sigma)$, i.e. the left action
\beq\label{eq:TacSi}
T\stac{L} \Sigma\to \Sigma,\qquad t\stac{L} x \mapsto t(x)
\eeq
is a coretraction in ${}_T\cM^\cD$.
\end{proposition}
\begin{proof} 
A retraction of the map \eqref{eq:TacSi} can be constructed 
in terms of the elements $\{\tildej_\ell\}\subset \widetilde{Q}$ and
$\{j_\ell\}\subset{_L\Hom^\Dd}(\Dd,\Sigma)$, satisfying
$\sum_\ell\tildej_\ell\newblackdiamondd j_\ell=\Cc$, 
as
$$
\sigma:\Sigma\to T\stac{L} \Sigma, \quad 
x\mapsto \sum_\ell
x^{[0]}\tildej_\ell({x^{[1]}}_{[0]})(-)\stac{L}j_\ell({x^{[1]}}_{[1]}) \equiv 
\sum_\ell
{x_{[0]}}^{[0]}\tildej_\ell({x_{[0]}}^{[1]})(-)\stac{L}j_\ell(x_{[1]}). 
$$
It is a well defined map by \leref{3.4} (2).
Being a composite of right $\cD$-colinear maps, it is
$\cD$-colinear. Its $T$-linearity follows by the fact that any $t\in T$ is
$\cC$-colinear and hence, in particular $\cD$-colinear and $A$-linear. That 
$\sigma$ is a retraction of the map \eqref{eq:TacSi} is a direct consequence of
\leref{jids} (2).
\end{proof}

\section{Weak and strong structure theorems}\label{sec:str.thm}

As an $L$-$\cC$ bicomodule $\Sigma$,
for a right coring extension $(\Dd:L)$ of $(\Cc:A)$, is in particular a right
$\cC$-comodule, it determines an adjunction of functors  
\beq\label{eq:ind}
-\stac{T} \Sigma\ :\ \cM_T\to \cM^\cC,
\eeq
from the category of right modules for the algebra $T=\mathrm{End}^\cC(\Sigma)$
to the category of right comodules for the coring $\cC$, as in the bottom row
of diagram \eqref{fig:EiMo}, and  
\beq\label{eq:coinv}
\Hom^\cC(\Sigma,-)\ :\ \cM^\cC\to \cM_T.
\eeq
In the present section we study the `descent theory' of pure coring extensions, 
i.e. investigate in what situations the functor \eqref{eq:coinv} is fully 
faithful or an equivalence with inverse \eqref{eq:ind}.
\begin{theorem}[Weak Structure Theorem]\label{thm:w.str.thm}\thlabel{w.str.thm}
Let the $L$-coring $\cD$ be a pure right extension of the $A$-coring
$\cC$. Take  
$\Sigma\in {}_L\cM^\cC$ and consider the associated Morita context 
$\widetilde{\mathbb M}(\Sigma)$ in  \equref{eq:contextmod}. If the map
$\newblackdiamond$ in \equref{blackdiamond} is surjective then the functor
\eqref{eq:coinv} is fully faithful. 
\end{theorem}
\begin{proof}
The property, that the functor \eqref{eq:coinv} is fully faithful,
is equivalent to the bijectivity of the counit of the adjunction,
\beq\label{eq:couM}
\varepsilon_M\ :\ \Hom^\cC(\Sigma,M)\stac{T} \Sigma \to M,\qquad
\varphi_M\stac{T} x\mapsto \varphi_M(x),
\eeq
for $T\colon=\mathrm{End}^\cC(\Sigma)$ and any right $\cC$-comodule $M$. Note
that the restriction of the map  $(M\sstac A \epsilon_\cC)\circ \can_M$ to  
$\Hom^\cC(\Sigma,M)\sstac{T} \Sigma$ is equal to $\varepsilon_M$.
Furthermore, by the right $\cC$-colinearity of $\varepsilon_M$, we have
 $\varrho^M\circ \varepsilon_M=
(\varepsilon_M\sstac A\cC)\circ (\mathrm{Hom}^\cC(\Sigma,M)\sstac
T\varrho^\Sigma)$, which equals the restriction of $\can_M$. By
\thref{surjectivity} (1), the 
surjectivity of $\newblackdiamond$ implies that $\Sigma$ is a Galois
  $\cC$-comodule. Taking the explicit form of $\can^{-1}$ in the
  proof of \thref{surjectivity} into account, we have 
$$
\can_M^{-1}\circ \varrho^M(m)=
\sum_\ell {m_{[0]}}^{[0]}\tildej_\ell({m_{[0]}}^{[1]})(-)\stac T
j_\ell(m_{[1]}),\qquad \textrm{for }m\in M,$$
which is an element of $\mathrm{Hom}^\cC(\Sigma,M)\sstac T \Sigma$, by
\leref{3.4} (2). In light of these observations,
the inverse of $\varepsilon_M$ can be constructed as 
$$
\can_M^{-1}\circ \varrho^M:M\to \Hom^\cC(\Sigma,M)\stac{T} \Sigma.
$$
\end{proof}

\begin{corollary}\label{cor:fgp.str.thm}
Let the $L$-coring $\cD$ be a pure right extension of the $A$-coring $\cC$ 
and $\Sigma$ an $L$-$\cC$ bicomodule. Consider the Morita contexts
$\widetilde{\mathbb M}(\Sigma)$, associated to
$\Sigma$ in \equref{eq:contextmod}, and ${\mathbb M}(\Sigma)$, 
associated to $\Sigma$ as a $\cC$-comodule in \eqref{eq:SiMor}.
Suppose that the map $\newblackdiamond$, given in \equref{blackdiamond}, is
surjective. 
Then the connecting map $\newblacktriangle$, given in \eqref{eq:F}, 
is surjective if and only if 
$\cC$ is a finitely generated and projective left $A$-module.
\end{corollary}
\begin{proof}
If the connecting map $\newblacktriangle$ in \eqref{eq:F} is
surjective then $\cC$ is a finitely generated projective left $A$-module by
Lemma \ref{lem:Afgp} (1).
Conversely, if $\cC$ is a finitely generated projective left $A$-module
then the 
connecting map \eqref{eq:F} is equal to the counit \eqref{eq:couM} for
the right $\cC$-comodule ${}^*\cC$ (cf. Remark \ref{rem:Qfgp}).
Then it is an isomorphism by Theorem \ref{thm:w.str.thm}.
\end{proof}

Recall from \thref{tensorff} that a sufficient condition for the functor
\eqref{eq:ind} to be fully faithful is the surjectivity of the map
\equref{defG}. Motivated by this result,
in what follows we look for conditions under which the map \equref{defG} is
surjective. 
\begin{proposition}\prlabel{diamondGsurj}
Let the $L$-coring $\Dd$ be a pure right extension of the $A$-coring $\Cc$ and
take $\Sigma\in{_L\Mm^\Cc}$. Consider the Morita contexts  
$\widetilde{\mathbb M}(\Sigma)$, associated to the $L$-$\cC$ bicomodule
$\Sigma$ in \equref{eq:contextmod}, and ${\mathbb M}(\Sigma)$,
associated to $\Sigma$ as a $\cC$-comodule in \eqref{eq:SiMor}.
In particular, denote $T=\mathrm{End}^\cC(\Sigma)$.
If the connecting map $\newdiamond$, given in \equref{diamond}, is surjective
and there exist elements 
$\{v_j\}\subset {}_L\mathrm{Hom}_L(\cD,T)$ and $\{d_j\}\subset \cD$, such 
that $\sum_j v_j(d_j)=1_T$, 
then also the connecting map $\newtriangle$, given in
\equref{defG}, is surjective (hence $\Sigma$ is a finitely generated and
projective right $A$-module by Lemma \ref{lem:Afgp} (2)).
\end{proposition}

\begin{proof}
By \leref{Qprop} (1), $\widetilde Q$ can be viewed as a $k$-submodule of
$Q\cong Q'$ in \eqref{eq:Q}. Hence (identifying $q\in \widetilde Q$ with the 
corresponding element of $Q$), it follows by the explicit forms of the maps
$\newdiamond$ and $\newtriangle$ that 
$(p\newdiamondd q)(d)=p(d)\newtriangled q$, for $p\in
{}_L\mathrm{Hom}^\cD(\cD,\Sigma)$, $q\in \widetilde{Q}$ and $d\in \cD$.
Let $\{h_i\}\subset {}_L\mathrm{Hom}^\cD(\cD,\Sigma)$ and
$\{\widetilde{h}_i\}\subset \widetilde{Q}$ be sets of morphisms such that
$\sum_i h_i\newdiamondd \widetilde{h}_i=\epsilon_\cD(-)1_T$. Then
\begin{eqnarray*}
1_T&=&\sum_j v_j(d_j)=\sum_jv_j(d_{j(1)})\epsilon_\cD(d_{j(2)})
=\sum_jv_j(d_{j(1)})\big(\sum_i(h_i\newdiamondd
\widetilde{h}_i)(d_{j(2)})\big)\\ 
&=&\sum_j
v_j(d_{j(1)})\big(\sum_i(h_i(d_{j(2)})\newtriangled\widetilde{h}_i\big) 
=\sum_i\big(\sum_j
v_j(d_{j(1)})(h_i(d_{j(2)}))\big)\newtriangled\widetilde{h}_i\\ 
&=&\sum_i\big(\sum_j(v_j h_i)(d_j)\big)\newtriangled\widetilde{h}_i,
\end{eqnarray*}
where the penultimate equality follows by the left $T$-linearity of
$\newtriangle$ and the last one follows by \equref{vp}.
Since $\newtriangle$ is $T$-$T$ bilinear, this proves its surjectivity.
\end{proof}

By standard Morita theory, \prref{diamondGsurj} implies the
following.  
\begin{corollary}\colabel{TsumSi}
Under the assumptions (and using the notations) of \prref{diamondGsurj}, 
$\Sigma$ is a generator of left $T$-modules. That is, $T$ is a direct summand
of a direct sum 
$\Sigma^z$, as a left $T$-module, where $z$ is the cardinality of the sets 
$\{h_i\}\subset {}_L\mathrm{Hom}^\cD(\cD,\Sigma)$ and
$\{\widetilde{h}_i\}\subset \widetilde{Q}$, such that $\sum_i
h_i\newdiamondd \widetilde{h}_i=\epsilon_\cD(-)1_T$.
\end{corollary}
\begin{remark}\label{rem:surj}
The assumption in \prref{diamondGsurj} about the existence of elements 
$\{v_j\}\subset {}_L\mathrm{Hom}_L(\cD,T)$ and $\{d_j\}\subset \cD$, such 
that $\sum_j v_j(d_j)=1_T$, holds in
various situations, studied in connection with cleft entwining structures in
\cite[Theorem 4.5]{CaVerWang:cleftentw} and \cite[Theorems 4.9 and
  4.10]{Abuh:cleftentw}.
\begin{zlist}
\item If the counit $\epsilon_\Dd$ of $\Dd$ is surjective, then 
there exists $d\in\Dd$ such that $\epsilon_\Dd(d)=1_L$. Putting $v:\Dd\to T$,
$d'\mapsto \epsilon_\Dd(d')1_T$, we have $v(d)=1_T$. 
\item If $\cD$ contains a grouplike element then it is mapped by
  $\epsilon_\cD$ to $1_L$, by definition. Hence $\epsilon_\cD$ is surjective,
  being $L$-$L$ bilinear, so the considerations in part (1) apply.
\item If $\cD$ is faithfully flat as a left, or as a right $L$-module then
  $\epsilon_\cD$ is surjective, since $\epsilon_\cD\sstac{L} \cD$ and
  $\cD\stac{L} \epsilon_\cD$ are epimorphisms, split by the coproduct
  $\Delta_\cD$. Hence this is an example of the situation in part (1) as
  well. 
\end{zlist}
\end{remark}
\begin{theorem}[Strong Structure Theorem]\label{thm:s.str.thm}
Let the $L$-coring $\cD$ be a pure right extension of the
$A$-coring $\cC$ and  
$\Sigma$ an $L$-$\cC$ bicomodule. Denote $T\colon =\mathrm{End}^\cC(\Sigma)$.
If the associated Morita context \equref{eq:contextmod} is strict and 
there exist elements
$\{v_j\}\subset {}_L\mathrm{Hom}_L(\cD,T)$ and $\{d_j\}\subset \cD$, such 
that $\sum_j v_j(d_j)=1_T$, 
then the functors \eqref{eq:ind} and \eqref{eq:coinv} are inverse
equivalences.  
\end{theorem}

\begin{proof}
This is an immediate consequence of \thref{w.str.thm}, \thref{tensorff}
and \prref{diamondGsurj}. 
\end{proof}

Note that under the hypothesis of Theorem \ref{thm:s.str.thm}, $\Sigma$ is a
finitely generated and projective right $A$-module (cf. Lemma \ref{lem:Afgp}
(2)), hence a Galois comodule in the sense of \cite{ElK:comcor}
(cf. \thref{surjectivity}). On the other hand, under the assumptions in
Theorem \ref{thm:s.str.thm}, $\Sigma$ is not necessarily flat as a left
$T$-module (equivalently, $\cC$ is not necessarily flat as a left
$A$-module). This way Theorem \ref{thm:s.str.thm} covers cases which are not
treated by the Galois Comodule Structure Theorem \cite[18.27]{BrzWis:cor}. This
will be clear by the observation in Section \ref{sec:ex}, that the Fundamental
Theorem of Hopf modules (for arbitrary Hopf algebras or Hopf algebroids) is a
particular instance of Theorem \ref{thm:s.str.thm}.

\section{Cleft bicomodules}\selabel{cleft}

The results in \seref{contextext} and Section \ref{sec:str.thm} allow for an
application to the main case of interest, when there exist `invertible'
elements in the Morita context, 
associated to a bicomodule of a pure coring extension, in the
following sense.  
Let $\cD$ be an $L$-coring which is a pure right extension of
an $A$-coring 
$\cC$. For an $L$-$\cC$ bicomodule $\Sigma$ consider the associated Morita
context  
$\widetilde{\mathbb M}(\Sigma)=\big(\ {}_L\mathrm{Hom}_L(\cD,T),$ 
${}^\cC\mathrm{End}^\cD(\cC)^{op}, {}_L\mathrm{Hom}^\cD(\cD,\Sigma),
\widetilde{Q},\newblackdiamond,\newdiamond\ \big)$ in
\equref{eq:contextmod}, where $T={\rm End}^\cC(\Sigma)$ and $\widetilde{Q}$ is
the $k$-module \eqref{eq:Q'}.
\begin{definition}\label{def:cleftcom}
An object $\Sigma$ of ${_L\Mm^\Cc}$ is called a {\em weak cleft} bicomodule
for the pure right coring extension $(\Dd:L)$ of $(\Cc:A)$ provided  
there exist elements $j\in {_L\Hom^\Dd}(\Dd,\Sigma)$ and $\tildej\in
\widetilde{Q}$ such that $\tildej\newblackdiamondd j=\cC$.  

An object $\Sigma$ of ${_L\Mm^\Cc}$ is called a {\em cleft} bicomodule
for the pure right coring extension $(\Dd:L)$ of $(\Cc:A)$
provided  
there exist elements $j\in {_L\Hom^\Dd}(\Dd,\Sigma)$ and $\tildej\in
\widetilde{Q}$ such that 
$\tildej\newblackdiamondd j=\cC$ and, in addition, $j\newdiamondd \tildej
=\epsilon_\cD(-)1_T$.
\end{definition}
Note that if $\Sigma\in {}_L\cM^\cC$ is a (weak) cleft bicomodule for a
pure right 
coring extension $(\cD:L)$ of $(\cC:A)$, with morphisms $j\in
{_L\Hom^\Dd}(\Dd,\Sigma)$ and $\tildej\in \widetilde{Q}$ as in Definition
\ref{def:cleftcom}, then, by \prref{contextmod}, the natural
transformation $\alpha_4(\tildej)$ in \prref{contextmod} is a (left)
inverse of $\alpha_3(j)$ there.

In standard Hopf Galois theory, cleft extensions are Galois extensions with
additional `normal basis property'. In order to derive a similar result for
pure coring extensions, we impose the following definition.
\begin{definition}
An $L$-$\cC$ bicomodule $\Sigma$ for a pure right coring
extension $(\cD:L)$ of 
$(\cC:A)$ is said to obey the  {\em weak normal basis property} if it is
isomorphic to a direct summand of $T\sstac{L}\Dd$, as a $T$-$\Dd$ bicomodule,
for $T=\mathrm{End}^\cC(\Sigma)$. 
 
An $L$-$\cC$ bicomodule $\Sigma$ for a pure right coring
extension $(\cD:L)$ of 
$(\cC:A)$ is said to obey the {\em normal basis property} if it is isomorphic
to $T\sstac{L} \cD$, as a $T$-$\Dd$ bicomodule.
\end{definition}
Note that the normal basis property of an $L$-$\cC$ bicomodule $\Sigma$, for a
pure right coring extension $(\cD:L)$ of $(\cC:A)$, implies the
isomorphism of the 
cotensor product $\Sigma\Box_\cD W$ to $T\sstac L W$, as a left $T\colon
=\End^\cC(\Sigma)$-module, for any left $\cD$-comodule $W$. This way,
properties (like projectivity or freeness) of the left $L$-module $W$ are
inherited by the associated left $T$-module $\Sigma\Box_\cD W$. 

It is immediately clear from the definition that the Morita context
\equref{eq:contextmod}, associated to a cleft bicomodule of a pure coring
extension, is strict. 
Hence the proof of \thref{surjectivity} and \reref{numbers} lead to the
following   
relation between cleft bicomodules and Galois comodules which satisfy the
normal basis property.  

\begin{corollary}
\colabel{th:3.2}\colabel{jJ}
(1) $\Sigma\in {}_L\cM^\cC$ is a weak cleft bicomodule for the pure right
coring extension  
  $(\Dd:L)$ of $(\Cc:A)$ if and only if $\Sigma$ is a Galois $\cC$-comodule 
  and satisfies the weak normal basis property; 

(2) $\Sigma\in {}_L\cM^\cC$ is a cleft bicomodule for the pure
  right coring extension $(\Dd:L)$ of $(\Cc:A)$ if and only if 
$\Sigma$ is a Galois $\cC$-comodule and satisfies the normal basis property.
\end{corollary}

\coref{TsumSi} has the following consequence.
\begin{corollary}
Let the $L$-coring $\cD$ be a pure right extension of the
$A$-coring $\cC$ and $\Sigma\in {}_L\cM^\cC$ a cleft bicomodule. Put $T\colon
=\mathrm{End}^\cC(\Sigma)$. If there exist elements 
$\{v_j\}\subset {}_L\mathrm{Hom}_L(\cD,T)$ and 
$\{d_j\}\subset \cD$ such that $\sum_j v_j(d_j)=1_T$ then $\Sigma$ contains
the left 
regular $T$-module as a direct summand. 
\end{corollary}
Theorems \ref{thm:w.str.thm} and \ref{thm:s.str.thm} imply the following
structure theorems.
\begin{corollary} \label{cor:cleft.str.thm}
For the functors \eqref{eq:ind} and \eqref{eq:coinv}, associated to a cleft
$L$-$\cC$ bicomodule $\Sigma$ of a pure right coring extension
$(\cD:L)$ of $(\cC:A)$, the following assertions hold.
\begin{zlist}
\item (Weak Structure Theorem) The functor \eqref{eq:coinv} is fully faithful.
\item (Strong Structure Theorem) The functors \eqref{eq:ind} and
  \eqref{eq:coinv} are inverse equivalences provided that there exist elements
$\{v_j\}\subset {}_L\mathrm{Hom}_L(\cD,T)$ and
  $\{d_j\}\subset \cD$, such that $\sum_j v_j(d_j)=1_T$, where the notation
  $T=\mathrm{End}^\cC(\Sigma)$ is used.
\end{zlist}
\end{corollary}

\section{Examples}\label{sec:ex}

The following proposition implies that most coring extensions occurring in
this section are pure without any further assumption.

\begin{proposition}\label{prop:split}
If an $L$-coring ${\mathcal D}$ is a right extension of an $A$-coring
${\mathcal C}$, such that the right $L$-action on ${\mathcal C}$ is induced by
an algebra map $L\to A$ from the right $A$-action on ${\mathcal C}$, then
${\mathcal D}$ is a pure right extension of ${\mathcal C}$.
\end{proposition}

\begin{proof}
Note that the same formula \eqref{eq:M}, if considered in the category 
of right $A$-modules, is a split equalizer in the
sense dual to \cite[p 110]{BarrWells:ttt}. The $A$-module splitting is given
by $M \ot_A 
{\mathcal C}\ot_A \varepsilon_{\mathcal C}: M\ot_A {\mathcal C} \ot_A
{\mathcal C}\to M\ot_A {\mathcal C}$, where $\varepsilon_{\mathcal C}$ is
the counit of ${\mathcal C}$. Under our hypothesis there is a forgetful
functor ${\mathcal M}_A\to {\mathcal M}_L$, that takes \eqref{eq:M} to  
a split equalizer in ${\mathcal M}_L$, for any right ${\mathcal C}$-comodule
$(M,\varrho)$.  
By a dual form of \cite[p 110 Proposition 2]{BarrWells:ttt}, split equalizers
are preserved by any functor, what proves the claim.  
\end{proof}

\subsection{Cleft entwining structures}
An {\em entwining structure} \cite{BrzMaj:coabund} consists of a
$k$-algebra $A$, a $k$-coalgebra $\cD$ and a $k$-linear map
$\psi:\cD\sstac{k} A \to A\sstac{k} \cD$, satisfying
\begin{eqnarray}
&&\psi\circ (\cD\stac{k}\mu_A)= (\mu_A\stac{k} \cD)\circ
(A\stac{k} \psi)\circ (\psi\stac{k}
  A)\label{eq:entwa}\\
&&\psi \circ (\cD\stac{k} 1_A)=1_A\stac{k} \cD \label{eq:entwb}\\
&&(A\stac{k} \Delta_\cD)\circ \psi= (\psi\stac{k} \cD)\circ
(\cD\stac{k} \psi)\circ (\Delta_\cD\stac{k} A)
  \label{eq:entwc} \\
&&(A\stac{k} \epsilon_\cD)\circ \psi=\epsilon_\cD\stac{k} A.
\label{eq:entwd}
\end{eqnarray}
The index notation $\psi(d\sstac{k} a)=a_\psi\sstac{k} d^\psi$ will be 
used, where implicit summation is understood. To an
entwining structure $(A,\cD,\psi)$ one can associate an $A$-coring
$\cC\colon =A\sstac{k} \cD$ as follows (\cite[Proposition
2.2]{Brz:str}). The left $A$-module structure is given by left
multiplication in the first tensorand and the right $A$-module
structure is given by $(a\sstac{k} d)a'=a a'_\psi\sstac{k} d^\psi$,
for $d\in \cD$ and $a,a'\in A$. The coproduct is
\begin{eqnarray*}
\Delta_\cC\colon = A\stac{k}\Delta_\cD:\cC\simeq A\stac{k}\cD
&\to&
\cC\stac{A}\cC\simeq A\stac{k}\cD\stac{k}\cD,\\
a\stac{k} d &\mapsto& (a\stac{k} d_{(1)})\stac{A} (1_A\stac{k}
d_{(2)})\simeq a\stac{k} d_{(1)} \stac{k} d_{(2)},
\end{eqnarray*}
and the counit is $\epsilon_\cC\colon = A\sstac{k} \epsilon_\cD$.
Clearly, $\cC$ is a $\cC$-$\cD$ bicomodule with the left regular
$\cC$-coaction $\Delta_\cC$ and right $\cD$-coaction
$\tau_\cC\colon = A\sstac{k} \Delta_\cD$. That is, $\cD$ is a right extension
of $\cC$, 
moreover, a pure coring extension by Proposition \ref{prop:split}.
This implies that any right $\cC$-comodule possesses a right
$\cD$-comodule structure. What is more, right 
$\cC$ comodules (also called {\em entwined modules}) are those
right $\cD$-comodules $M$ that are right $A$-modules as well and
the compatibility condition 
\beq\label{eq:entwmod}
(ma)_{[0]}\stac{k}(ma)_{[1]}=m_{[0]}a_\psi\stac{k}{m_{[1]}}^\psi
\eeq 
holds, for any $m\in M$ and $a\in A$.

An entwining structure $(A,\cD,\psi)$ has been termed {\em cleft}
in \cite[Definition 4.6]{Abuh:cleftentw} if $A$ (with the right
regular $A$-module structure) is an entwined module and there
exists a convolution invertible right $\cD$-colinear map
$\lambda:\cD\to A$.
(In the paper \cite{CaVerWang:cleftentw} a cleft entwining
structure is meant in the more restrictive sense that in addition
$\cD$ possesses a grouplike element $x$ and the right
$\cD$-coaction in $A$ is of the form $a\mapsto \psi(x\sstac{k} a)$.
Note that in this case $\End^\cC(A)\simeq\Hom^\cD(k,A)$, i.e. the $\cC$, 
and $\cD$-coinvariants of $A$ coincide.)
\begin{proposition} \label{prop:cleftentw}
An entwining structure $(A,\cD,\psi)$ is cleft if and only if $A$
(with the right regular $A$-module structure) is a cleft
bicomodule for the pure coring extension $\cD$ of $\cC\colon =
A\sstac{k}\cD$. 
\end{proposition}
\begin{proof}
Let us assume first that $(A,\cD,\psi)$ is a cleft entwining
structure. In this case $A$ is an entwined module, i.e. a
right $\cC$-comodule, by assumption. Let $\lambda:\cD\to A$ be a right
$\cD$-colinear map, with convolution inverse ${\bar \lambda}$. Put
$j\colon = \lambda$ and ${\tildej}:\cC\to A$, $a\otimes_k
d\mapsto a{\bar \lambda}(d)$. We need to prove that ${\tildej}$
is an element of the appropriate bimodule \eqref{eq:Q'} in the Morita 
context ${\widetilde{\mathbb M}}(A)$, associated to $A$ as in
\equref{eq:contextmod}, that is, of
$$
\widetilde{Q}\simeq \{q\in {}_A {\rm Hom}(\cC,A)\ |\
\forall d\in \cD,a\in A\quad 
\psi\big(d_{(1)}\stac{k} 
q(1_A\stac{k} d_{(2)})a\big)=q(1_A\stac{k}d)a_{[0]}\stac{k}
a_{[1]}\},
$$
(which is equal to ${}^\cC\Hom(\cC,A)$, cf. Remark \ref{rem:Sifgp}).
Note that, by \cite[Lemma 4.7 1]{Abuh:cleftentw}, the identity
\beq\label{eq:lam} 
{\bar \lambda}(d) 1_{A[0]}\stac{k} 1_{A[1]}={\bar \lambda}(d_{(2)})_\psi
  \stac{k}{d_{(1)}}^\psi
\eeq 
holds true, for any $d\in \cD$. Using the assumption that $A$
is an entwined module (in the second equality), \eqref{eq:lam} (in
the third one), and property \eqref{eq:entwa} of entwining
structures (in the fourth one), one checks that, for $d\in \cD$
and $a\in A$,
\begin{eqnarray}\label{eq:cleftd}
{\tildej}(1_A\stac{k} d)a_{[0]}\stac{k} a_{[1]} &=&{\bar
\lambda}(d) a_{[0]}\stac{k} a_{[1]} ={\bar \lambda}(d)
1_{A[0]}a_\psi \stac{k} {1_{A[1]}}^\psi ={\bar
\lambda}(d_{(2)})_{\psi'} a_\psi \stac{k} {d_{(1)}}^{\psi'\psi}
\nonumber\\
&=& \psi\big(d_{(1)}\stac{k} {\bar \lambda}(d_{(2)})a\big)
=\psi\big(d_{(1)}\stac{k} {\tildej}(1_A\stac{k} d_{(2)})a\big),
\end{eqnarray}
that is, $\tildej\in \widetilde{Q}$.
By the assumption that ${\bar \lambda}$ is left convolution
inverse of $\lambda$, for $a\in A$ and $d\in \cD$,
\begin{eqnarray}\label{eq:centw} 
({\tildej}\newblackdiamondd j)(a\stac{k} d)
&=&a{\bar \lambda}(d_{(1)})\lambda(d_{(2)})_{[0]}\stac{k}
\lambda(d_{(2)})_{[1]}\\
&=& a{\bar \lambda}(d_{(1)})\lambda(d_{(2)})\stac{k}
d_{(3)}=a\epsilon_\cD(d_{(1)}) \stac{k} d_{(2)}=a\stac{k} d,\nonumber
\end{eqnarray}
where the second equality follows by the colinearity of $\lambda$.
Similarly, since ${\bar \lambda}$ is also right convolution
inverse of $\lambda$, for $d\in \cD$,
$$
(j\newdiamond {\tildej})(d)= \lambda(d)_{[0]}{\bar
\lambda}\big(\lambda(d)_{[1]}\big) =\lambda(d_{(1)}){\bar
\lambda}(d_{(2)}) =\epsilon_\cD(d)1_A.
$$
This proves that $A$ is a cleft bicomodule, as stated.

Conversely, assume that $A$ is a cleft bicomodule for the
pure coring extension $\cD$ of $\cC$. Then it is, in
particular, an entwined module. Let $j\in \mathrm{Hom}^\cD(\cD,A)$ and
${\tildej}\in \widetilde{Q}$ be elements of the 
bimodules in the Morita context 
${\widetilde{\mathbb M}}(A)$, associated to $A$ as in 
\equref{eq:contextmod}, such that ${\tildej}\newblackdiamondd j=\cC$ and
$j\newdiamondd 
{\tildej}=1_A\epsilon_\cD(-)$. Then $\lambda\colon =j:\cD\to A$
is right $\cD$-colinear and ${\bar \lambda}:d\mapsto {\tildej}(1_A\sstac{k} d)$
is its convolution inverse. 
\end{proof}

\subsection{Cleft extensions of algebras by a coalgebra}
Let $\cD$ be a coalgebra over $k$ and $A$ a $k$-algebra and a right
$\cD$-comodule. In \cite{Brz:cogal} $A$ has been termed a {\em
$\cD$-cleft} extension of the subalgebra
$$
T\colon =\{\ t\in A\ |\ \forall a\in A\quad (ta)_{[0]}\stac{k}
(ta)_{[1]} = t a_{[0]}\stac{k} a_{[1]}\ \},
$$
if it is a $\cD$-Galois extension, i.e. the canonical map
\beq\label{eq:coacan} 
A\stac{T} A\to A\stac{k} \cD, \qquad a\stac{T} a'\mapsto aa'_{[0]}\stac{k}
a'_{[1]}  
\eeq 
is bijective,
and there exists a convolution invertible right $\cD$-colinear map
$\lambda:\cD\to A$.

Recall that for any $\cD$-Galois extension $A$ of $T$ there exists
a (unique) entwining structure $(A,\cD,\psi)$ such that $A$ is an
entwined module (cf. \cite[34.6]{BrzWis:cor}). On the other hand,
the canonical map \eqref{eq:coacan} is bijective for any cleft
entwining structure $(A,\cD,\psi)$ by \cite[Proposition 4.8
1]{Abuh:cleftentw}. This means that cleft extensions of algebras by a 
coalgebra are in one-to-one correspondence with cleft entwining 
structures. Combining this observation with Proposition
\ref{prop:cleftentw}, we conclude that $A$ is a $\cD$-cleft
extension of $T$ if and only if $A$ is a cleft bicomodule for the
pure coring extension $\cD$ of $\cC\colon = A\sstac{k} \cD$.

\subsection{Cleft extensions of algebras by a Hopf algebra}
Let $\cD$ be a Hopf algebra over $k$ and $A$ a right comodule
algebra. The algebra $A$ and the coalgebra underlying $\cD$ are
entwined by the map
$$
\psi:\cD\stac{k} A\to A\stac{k} \cD,\qquad d\stac{k} a\mapsto
a_{[0]}\stac{k} d a_{[1]}.
$$
Since $1_\cD$ is a grouplike element in $\cD$, $1_A\sstac{k} 1_\cD$
is a grouplike element in the $A$-coring $\cC\colon = A\sstac{k}
\cD$, associated to the entwining structure $(A,\cD,\psi)$. Hence
$A$ is an entwined module.

$A$ is called a {\em $\cD$-cleft} extension of its
$\cD$-coinvariant subalgebra if and only if there exists a
convolution invertible right $\cD$-colinear map $\lambda:\cD\to A$
(see e.g. \cite[Definition 7.2.1]{Montg:hopfa}), i.e. if and only
if $(A,\cD,\psi)$ is a cleft entwining structure. 
(Note that this way a cleft extension of algebras by a Hopf algebra 
is a cleft  
extension by the underlying coalgebra.) By Proposition
\ref{prop:cleftentw}, this is equivalent to $A$ being a cleft
bicomodule for the pure coring extension $\cD$ of $\cC\colon =
A\sstac{k} \cD$. 

\subsection{Cleft weak entwining structures}
A {\em weak entwining structure} \cite{CaDeGr:wentw} consists of a
$k$-algebra $A$, a $k$-coalgebra $\cD$ and a $k$-linear map
$\psi:\cD\sstac{k} A \to A\sstac{k} \cD$, such that the
compatibility conditions \eqref{eq:entwa} and \eqref{eq:entwc}
hold true, while \eqref{eq:entwb} and \eqref{eq:entwd} are
replaced by
\begin{eqnarray}
&&\psi\circ(\cD\stac{k} 1_A)=(e\stac{k}\cD)\circ \Delta_\cD,\qquad
  \textrm{and}\label{eq:wentwb}\\
&&(A\stac{k} \epsilon_\cD)\circ \psi=\mu_A\circ (e\stac{k}A),
\label{eq:wentwd}
\end{eqnarray}
respectively, where $e\colon =(A\sstac{k} \epsilon_\cD)\circ
\psi\circ (\cD\sstac{k} 1_A):\cD\to A$.

To a weak entwining structure $(A,\cD,\psi)$ one can associate an
$A$-coring $\cC\colon =\{\ a 1_{A\psi}\sstac{k} d^\psi\
\}_{a\sstac{k} d\in A\sstac{k}\cD}$ (cf. \cite[Proposition 2.3]{Brz:str} 
or \cite[37.4]{BrzWis:cor}).
The left $A$-module structure is given by left multiplication in
the first tensorand, and the right $A$-module structure is given by $(a
1_{A\psi}\sstac{k} d^\psi)a'=aa'_\psi \sstac{k} d^\psi$. The
coproduct is given by the restriction of $A\sstac{k}\Delta_\cD$,
i.e. by
\begin{eqnarray*}
\Delta_\cC:\quad\cC\quad&\to& \cC\stac{A}\cC,\\
a 1_{A\psi}\stac{k} d^\psi &\mapsto& (a 1_{A\psi}\stac{k}
{d^\psi}_{(1)})\stac{A} (1_A\stac{k} {d^\psi}_{(2)})= (a
1_{A\psi}\stac{k} {d_{(1)}}^\psi)\stac{A}( 1_{A\psi'}\stac{k}
{d_{(2)}}^{\psi'}).
\end{eqnarray*}
The counit is given by the restriction of $A\sstac{k}
\epsilon_\cD$, i.e. by
$$
\epsilon_\cC:\cC\to A,\qquad a 1_{A\psi}\stac{k} d^\psi\mapsto a
1_{A\psi}\epsilon_\cD(d^\psi)= ae(d).
$$
$\cC$ is a $\cC$-$\cD$ bicomodule with the left regular
$\cC$-coaction $\Delta_\cC$ and right $\cD$-coaction, given by the
restriction of $A\sstac{k}\Delta_\cD$, i.e.
\begin{equation}\label{eq:tauC}
\tau_\cC:\cC\to \cC\stac{k}\cD,\quad a 1_{A\psi}\stac{k}
d^\psi\mapsto (a 1_{A\psi}\stac{k} {d^\psi}_{(1)}) \stac{k}
{d^\psi}_{(2)}= (a 1_{A\psi}\stac{k} {d_{(1)}}^\psi)\stac{k}
d_{(2)},
\end{equation}
where the equality of the two forms of $\tau_\cC$ follows by \eqref{eq:wentwb}
and the coassociativity of $\Delta_\cD$.
That is, $\cD$ is a right extension of $\cC$. 
This coring extension is pure by Proposition \ref{prop:split}.
Right
$\cC$-comodules are called {\em weak entwined modules} and they
can be characterized as right $\cD$-comodules $M$, that are right
$A$-modules as well such that the compatibility condition
\eqref{eq:entwmod} holds true.

By \cite[Definition 1.9]{AAFVGRRR:wcleft}, a weak entwining
structure $(A,\cD,\psi)$ is {\em cleft} if $A$ (with the right
regular $A$-module structure) is a weak entwined module and there
exists a right $\cD$-colinear map $\lambda:\cD\to A$ and a
$k$-linear map ${\bar
  \lambda}:\cD\to A$, satisfying \eqref{eq:lam} and
\beq 
1_{A\psi}{\bar \lambda}(d^\psi)={\bar \lambda}(d)\qquad
\textrm{and}\qquad {\bar \lambda}(d_{(1)})\lambda(d_{(2)})=e(d),
\qquad \textrm{for }d\in \cD. 
\label{eq:wlinv} 
\eeq 
(The first condition in \eqref{eq:wlinv} can be read as a convenient
normalization. Indeed, if there exists ${\bar \lambda}\in {\rm
Hom}_k(\cD,A)$, satisfying \eqref{eq:lam} and the second condition in
\eqref{eq:wlinv}, then it can be replaced by the (non-zero) map
$d\mapsto 1_{A\psi}{\bar \lambda}(d^\psi)$.)
\begin{proposition}\label{prop:wcleft}
A weak entwining structure $(A,\cD,\psi)$ is cleft if and only if
$A$ (with the right regular $A$-module structure) is weak cleft bicomodule for
the pure
right coring extension $\cD$ of the $A$-coring $\cC$, associated to the
weak entwining structure $(A,\cD,\psi)$.
\end{proposition}
\begin{proof}
Let us assume first that $(A,\cD,\psi)$ is a cleft weak entwining
structure. We construct elements $j$ and ${\tildej}$ in the
bimodules of the Morita context \equref{eq:contextmod}, associated to
$A$, such that ${\tildej}\newblackdiamondd j=\cC$. Put $j\colon =
\lambda$ and
$$
{\tildej}:\cC\to A,\qquad  a 1_{A\psi}\stac{k} d^\psi\mapsto
a1_{A\psi}{\bar \lambda}(d^\psi)= a{\bar\lambda}(d),
$$
where $\lambda:\cD\to A$ is a right $\cD$-colinear map and ${\bar
\lambda}:\cD\to A$ is a $k$-linear map, satisfying \eqref{eq:lam}
and \eqref{eq:wlinv}. Analogously to \eqref{eq:cleftd} and
\eqref{eq:centw}, assumption \eqref{eq:lam} implies that ${\tildej}$ is 
left $\cC$-colinear, i.e. an element of 
\begin{eqnarray*}
\widetilde{Q}\simeq \{q\in {}_A {\rm Hom}(\cC,A) &|& \forall d\in \cD,a\in A
\nonumber\\
&&\psi\big(d_{(1)}\stac k q(1_{A\psi}\stac k
{d_{(2)}}^\psi)a\big)=q(1_{A\psi}\stac kd^\psi)a_{[0]}\stac k
a_{[1]} \ \},
\end{eqnarray*}
and \eqref{eq:wlinv} implies ${\tildej}\newblackdiamondd j=\cC$.

Conversely, assume that $A$ is a weak cleft bicomodule, i.e. there exist
elements $j\in \mathrm{Hom}^\cD(\cD,A)$ and ${\tildej}\in \widetilde{Q}$ in
the bimodules of the Morita 
context $\widetilde{\mathbb M}(A)$, associated to
$A$ in \equref{eq:contextmod}, such that 
${\tildej}\newblackdiamondd j=\cC$. The map $\lambda\colon = j:\cD\to A$ is 
right $\cD$-colinear. Together with the map ${\bar \lambda}:d\mapsto
{\tildej}(1_{A\psi}\sstac{k} d^\psi)$, for $d\in \cD$, they
satisfy \eqref{eq:wlinv} and \eqref{eq:lam}. Indeed, the first
condition in \eqref{eq:wlinv} follows by the left $A$-linearity of
${\tildej}$ and \eqref{eq:entwa}. The second one follows by
${\tildej}\newblackdiamondd j=\cC$ as, for $d\in \cD$,
\begin{eqnarray*}
e(d)&=&1_{A\psi} \epsilon_\cD(d^\psi) =
(A\stac k\epsilon_\cD)
\big((\tildej\newblackdiamondd j)(1_{A\psi}\stac k d^\psi)\big)\\
&=&{\tildej}(1_{A\psi}\stac{k} {d_{(1)}}^\psi)j( {d}_{(2)})
={\bar \lambda}(d_{(1)})\lambda (d_{(2)}),
\end{eqnarray*}
where the third equality follows by the forms \equref{blackdiamond} of
the map $\newblackdiamond$ and \eqref{eq:tauC} of the $\cD$-coaction in
$\cC$. 
Condition \eqref{eq:lam} is easily seen to follow by the assumption that 
${\tildej}$ is an element of
the bimodule $\widetilde{Q}$. 
\end{proof}
Note that the first condition in \eqref{eq:wlinv} and \eqref{eq:lam}, imposed
on a weak entwining 
structure in \cite{AAFVGRRR:wcleft}, gain an explanation by
Proposition \ref{prop:wcleft}. They mean that ${\bar \lambda}\in
\mathrm{Hom}_k(\cD,A)$ 
corresponds to an element of $\widetilde{Q}\subseteq
{}_A\mathrm{Hom}(\cC,A)$ in Proposition \ref{prop:wcleft}, via the isomorphism 
$$
{}_A{\rm Hom}(\cC,A)\simeq\{\ \nu\in{\rm
Hom}_k(\cD,A)\ |\ \forall d\in
\cD\quad 1_{A\psi}\nu(d^\psi)=\nu(d)\ \}.
$$

\subsection{Cleft extensions by partial group actions}
Extending the definition of (idempotent) partial actions of finite
groups on commutative algebras in \cite{DoEx:partac} and
\cite{DoFePa:partgal}, Caenepeel and De Groot introduced in
\cite{CaDeGr:partgal} idempotent partial actions of finite groups
$G$ on arbitrary algebras $A$, as follows. An idempotent partial
$G$-action on $A$ consists of a collection
$\{e_\sigma\}_{\sigma\in G}$  of central idempotents in $A$ and a
collection $\{\alpha_\sigma:Ae_{\sigma^{-1}}\to
Ae_\sigma\}_{\sigma\in G}$ of isomorphisms of ideals, satisfying
the conditions
\begin{eqnarray*}
&& A_1=A\quad \textrm{and}\quad \alpha_1=A,\qquad \qquad \qquad
\,\quad
\textrm{for the unit element }1\textrm{ of }G,\textrm{ and}\\
&&\alpha_\sigma\big(\alpha_\tau(a
e_{\tau^{-1}})e_{\sigma^{-1}}\big)= \alpha_{\sigma\tau}(a
e_{\tau^{-1}\sigma^{-1}})e_\sigma, \quad \textrm{for
}\sigma,\tau\in G,\ a\in A.
\end{eqnarray*}
They constructed an $A$-coring for such a partial action,
$\{e_\sigma,\alpha_\sigma\}_{\sigma\in G}$ of $G$ on $A$, as a
$k$-module $\cC\colon =\oplus_{\sigma\in G} Ae_\sigma$ with
$A$-$A$ bimodule structure
$$
a_1(a\nu_\sigma)a_2=a_1 a \alpha_\sigma(a_2
e_{\sigma^{-1}})\nu_\sigma,
$$
for $a_1,a_2,a\in A$, and elements $\nu_\sigma$ of $\cC$, taking
the value $e_\sigma$ in the component $\sigma$ and $0$ everywhere
else, for ${\sigma\in G}$. The coproduct and the counit are
inherited from the coalgebra $k(G)$, the $k$-dual of the group
algebra. Explicitly,
$$
\Delta_\cC(a\nu_\sigma)=\sum_{\tau\in G} a\nu_\tau \stac{A}
\nu_{\tau^{-1}\sigma} \qquad \textrm{and}\qquad
\epsilon_\cC(a\nu_\sigma)=a\delta_{\sigma,1}, \qquad \textrm{for
}a\nu_\sigma\in \cC.
$$
Note that the coalgebra ($k$-coring) $k(G)$ is a pure right
extension of 
the $A$-coring $\cC$. That is, there exists a left $\cC$-colinear
right $k(G)$-coaction in $\cC$,
$$
\tau_\cC: \cC\to \cC\stac{k} k(G),\qquad a\nu_\sigma \mapsto
\sum_{\tau\in G} a \nu_\tau\stac{k}\ u_{\tau^{-1}\sigma},
$$
where $\{u_{\sigma}\}_{\sigma \in G}$ is the $k$-basis for $k(G)$,
dual to the basis $\{\sigma\}_{\sigma\in G}$ of the group algebra.
Since $\cC$ possesses a grouplike element, \beq\label{eq:grlike}
\sum_{\sigma\in G}\nu_\sigma, \eeq (cf. \cite[Lemma
2.3]{CaDeGr:partgal}), $A$ possesses a right $\cC$-comodule
structure (and hence a right $k(G)$-comodule structure).

By a
plausible definition we call $A$ a {\em cleft extension} of its
$G$-invariant subalgebra $\{\ a\in A \ |\ \forall \sigma\in G\quad 
\alpha_\sigma (ae_{\sigma^{-1}})=ae_\sigma\ \}$ if
there exists a convolution invertible right colinear map from the
right regular $k(G)$-comodule to $A$.
\begin{proposition}
Let $G$ be a finite group with an idempotent partial action
$\{e_\sigma,\alpha_\sigma\}_{\sigma\in G}$ on an algebra $A$. Let
$\cC$ be the associated $A$-coring. Then $A$ is a cleft extension
of its $G$-invariant subalgebra if and only if $A$ is a cleft
bicomodule for the pure coring extension $k(G)$ of $\cC$.
\end{proposition}
\begin{proof}
The proof is surprisingly similar to that of Proposition
\ref{prop:cleftentw}.

Assume first that $A$ is a cleft bicomodule, that is, there
exist elements ${\tildej}\in\widetilde{Q}$ and $j\in
\mathrm{Hom}^{k(G)}(k(G),A)$ in the bimodules of the Morita  
context $\widetilde{\mathbb M}(A)$, associated to
$A$ as in \equref{eq:contextmod}, such that
${\tildej}\newblackdiamondd j=\cC$ and $j\newdiamondd {\tildej}= 
1_A\epsilon_{k(G)} (-)$. Put $\lambda\colon = j:k(G)\to A$. It is
right colinear. We claim that it is also convolution invertible.

Using the notations, introduced earlier in this section, the
$k(G)$-coaction in $A$, determined by the grouplike element
\eqref{eq:grlike}, comes out as
$$
\tau_A:A\to A\stac{k} k(G),\qquad a\mapsto \sum_{\sigma\in G}
\alpha_\sigma (a e_{\sigma^{-1}})\stac{k}\ u_\sigma.
$$
Then, since $k(G)$ is a free $k$-module of finite rank, the
colinearity of $\lambda$ means that \beq\label{eq:jcolin}
\alpha_\tau\big(\lambda(u_\sigma)
e_{\tau^{-1}}\big)=\lambda(u_{\sigma\tau^{-1}}), \qquad
\textrm{for }\sigma,\tau\in G. \eeq Condition \eqref{eq:jcolin}
implies, in particular, that
$\lambda(u_\sigma)e_\tau=\lambda(u_\sigma)$, for any
$\sigma,\tau\in G$. (Hence there exist no non-trivial right
$k(G)$-comodule maps $k(G)\to A$ if the ideals
$\{Ae_\sigma\}_{\sigma \in G}$ have no non-trivial intersection.)

Now put ${\bar \lambda}(u_\sigma)\colon = {\tildej}(\nu_\sigma)$,
for $\sigma\in G$. By ${\tildej}\newblackdiamondd j=\cC$,
${\bar \lambda}$ is left convolution inverse of $\lambda$.
Similarly, it follows by $j\newdiamondd {\tildej}=1_A\epsilon_{k(G)}
(-)$ that
$$
\sum_{\tau\in G} \alpha_\tau\big( \lambda(u_\sigma)
e_{\tau^{-1}}\big) {\bar \lambda}(u_\tau)=\delta_{\sigma,1}\, 1_A\qquad
\textrm{for }\sigma \in G.
$$
Using the colinearity of $\lambda$, i.e. the identity
\eqref{eq:jcolin}, we conclude that ${\bar \lambda}$ is also right
convolution inverse of $\lambda$.

Conversely, assume that there exists a right $k(G)$-comodule map
$\lambda:k(G)\to A$ with convolution inverse ${\bar \lambda}$. We
construct elements $j\in \mathrm{Hom}^{k(G)}(k(G),A)$ and ${\tildej}\in
\widetilde{Q}$ such that
${\tildej}\newblackdiamondd j=\cC$ and $j\newdiamondd {\tildej}=
1_A\epsilon_{k(G)} (-)$. Put $ j\colon =\lambda:k(G)\to A$. Since
${\bar \lambda}$ is convolution inverse of a right $k(G)$-comodule
map $\lambda$, its range is in the intersection of the ideals
$\{Ae_\sigma\}_{\sigma\in G}$. Hence we can put
$$
{\tildej}: \cC\to A,\qquad a\nu_\sigma \mapsto a{\bar
\lambda}(u_\sigma).
$$
The conditions ${\tildej}\newblackdiamondd j=\cC$ and $j\newdiamondd
{\tildej}=1_A\epsilon_{k(G)} (-)$ follow easily by the
assumptions that ${\bar \lambda}$ is left, and right convolution
inverse of $\lambda$, respectively, and the colinearity condition
\eqref{eq:jcolin}. Furthermore, using that ${\bar \lambda}$ is
left convolution inverse of $\lambda$ (in the second equality),
the colinearity condition \eqref{eq:jcolin} (in the third one) and
the assumption that ${\bar \lambda}$ is left convolution inverse
of $\lambda$ (in the last one), we
deduce that 
\begin{eqnarray*}
\alpha_\tau\big( {\bar
\lambda}(u_{\tau^{-1}\sigma})ae_{\tau^{-1}}\big) &=&
\sum_{\omega\in G} \delta_{\omega,\tau}
\alpha_\tau\big( {\bar \lambda}(u_{\omega^{-1}\sigma})ae_{\tau^{-1}}\big)\\
&=& \sum_{\omega,\mu\in G} {\bar
\lambda}(u_\mu)\lambda(u_{\mu^{-1}\omega\tau^{-1}})
\alpha_\tau\big( {\bar \lambda}(u_{\omega^{-1}\sigma})ae_{\tau^{-1}}\big)\\
&=& \sum_{\omega,\mu\in G} {\bar \lambda}(u_\mu) \alpha_\tau\big(
\lambda(u_{\mu^{-1}\omega})e_{\tau^{-1}}\big)
\alpha_\tau\big( {\bar \lambda}(u_{\omega^{-1}\sigma})ae_{\tau^{-1}}\big)\\
&=& \sum_{\omega,\mu\in G} {\bar \lambda}(u_\mu) \alpha_\tau\big(
\lambda(u_{\mu^{-1}\omega}) {\bar
\lambda}(u_{\omega^{-1}\sigma})ae_{\tau^{-1}}\big)
={\bar \lambda}(u_\sigma) \alpha_\tau(ae_{\tau^{-1}}),
\end{eqnarray*}
for $a\in A$ and $\sigma,\tau\in G$. Using the forms of the coproduct
$\Delta_\cC$ in $\cC$ and the $\cC$-coaction (determined by the grouplike
element \eqref{eq:grlike}) in $A$, it is straightforward to 
check that this is equivalent to the property that ${\tildej}$ is
an element of the bimodule $\widetilde{Q}$, associated to $A$ as in
\eqref{eq:Q'}. 
\end{proof}

\subsection{Cleft entwining structures over arbitrary base}\label{sec:Lcleft}
An {\em entwining structure over an algebra $L$} consists of an
$L$-ring $A$, an $L$-coring $\cD$ and an $L$-$L$ bilinear map
$\psi:\cD\sstac{L} A\to A\sstac{L}\cD$, satisfying conditions
(\ref{eq:entwa}-\ref{eq:entwd}), with the only modification that
$k$-module tensor products are replaced by $L$-module tensor
products. Just as in the case of commutative base rings,
$\cC\colon = A\sstac{L} \cD$ possesses an $A$-coring structure (cf.
\cite[Example 4.5]{Bohm:int}) such that $\cD$ is a 
pure right extension of $\cC$.

Recall that, for an $L$-ring $A$ and an $L$-coring $\cD$, the set of 
bimodule maps ${}_L\Hom_L(\cD,A)$ is an algebra with the convolution 
product $(fg)(d)=f(d_{(1)})g(d_{(2)})$ and unit $\epsilon_\cD(-)1_A$. 
In complete analogy with
Proposition \ref{prop:cleftentw} one proves the following.
\begin{proposition}\label{prop:cleftLentw}
Let $(A,\cD,\psi)$ be an entwining structure over an algebra $L$
and let $\cC\colon =A\sstac{L} \cD$ be the associated $A$-coring.
$A$ (with the right regular $A$-module structure) is a cleft
$L$-$\cC$ bicomodule for the pure coring extension $\cD$ of
$\cC$ if 
and only if the following assertions hold.
\begin{blist}
\item $A$ (with the right regular $A$-module structure) is an
entwined module, i.e. it is a right $\cD$-comodule such that the
  compatibility condition
$$
(aa')_{[0]}\stac{L} (aa')_{[1]}=a_{[0]}a'_\psi\stac{L}
{a_{[1]}}^\psi
$$
holds true, for all $a,a'\in A$; 
\item The right $\cD$-coaction $a\mapsto a_{[0]}\sstac{L} a_{[1]}$ 
in $A$ is left $L$-linear;
\item There exists a convolution invertible morphism $\lambda\in
  {}_L\Hom^\cD(\cD,A)\subseteq {}_L\mathrm{Hom}_L(\cD,A)$.
\end{blist}
\end{proposition}
If conditions (a)-(c) in Proposition \ref{prop:cleftLentw} hold,
we call the $L$-entwining structure $(A,\cD,\psi)$ {\em cleft}.

Note that if the coring $\cD$ possesses a grouplike element $x$,
then the left $L$-linearity of the $\cD$-coaction $a\mapsto
\psi(x\sstac{L} a)$ in $A$ is equivalent to the element $x$ to be
central in the $L$-$L$ bimodule $\cD$.

\subsection{Cleft extensions of algebras by a Hopf algebroid}
A {\em Hopf algebroid} ${\mathcal H}$  consists of a left
bialgebroid structure ${\mathcal H}_L$, over a base algebra $L$,
and a right bialgebroid structure ${\mathcal H}_R$, over a base
algebra $R$, on the same total algebra $H$, and a $k$-linear map
$S:H\to H$, called the antipode, relating the two bialgebroid
structures \cite{BohmSzl:hgdax},\cite{Bohm:hgdint}. For a Hopf 
algebroid ${\mathcal H}$, we denote by $\gamma_L$ and 
$\pi_L$ (resp. $\gamma_R$ and $\pi_R$) the coproduct and the counit 
of the constituent bialgebroid ${\mathcal H}_L$ (resp. ${\mathcal H}_R$).

In a Hopf algebroid, the coproduct $\gamma_L:H\to H\ot_L H$ is a right
${\mathcal H}_R$-comodule map and $\gamma_R:H\to H\ot_R H$ is a right
${\mathcal H}_L$-comodule map. That is, (the coring underlying) ${\mathcal
  H}_R$ is a right extension of (the coring underlying) ${\mathcal H}_L$ and
(the coring underlying) ${\mathcal H}_L$ is a right extension of (the coring
underlying) ${\mathcal H}_R$. 
Hopf algebroids, in which both of these coring extensions are pure,
are called {\em pure Hopf algebroids}.
We are not aware of any examples of Hopf algebroids that are not pure. 

For a pure Hopf algebroid ${\mathcal H}$, the 
categories of right ${\mathcal H}_L$-comodules and of ${\mathcal
  H}_R$-comodules are isomorphic monoidal categories such that the forgetful
functor to the bimodule category ${}_R\cM_R$ is strict monoidal, 
cf. \cite[\emph{Corrigendum}, Theorems 4 and 6]{BohmBrz:cleft}.
Right ${\mathcal H}_R$-comodule algebras are defined as monoids in the
category of right ${\mathcal H}_R$-comodules, hence they are in
particular $R$-rings \cite{Scha:bianc}. A right ${\mathcal
H}_R$-comodule algebra $A$ determines an entwining structure over
$R$ with $R$-ring $A$ and $R$-coring $(H,\gamma_R,\pi_R)$,
underlying the bialgebroid ${\mathcal H}_R$ (cf.
\cite[(3.17)]{Bohm:gal}). Hence there exists a corresponding $A$-coring
$\cC\colon = A\sstac{R} H$, with coproduct $A\sstac{R} \gamma_R$,
inherited from ${\mathcal H}_R$. By right ${\mathcal H}_L$-colinearity of
$\gamma_R$, $\cC$ possesses a
$\cC$-${\mathcal H}_L$ bicomodule structure with the left regular
$\cC$-coaction and right
${\mathcal H}_L$-coaction $A\sstac{R} \gamma_L$. That is, the
$L$-coring $(H,\gamma_L,\pi_L)$, underlying the bialgebroid
${\mathcal H}_L$, is a right extension of the $A$-coring 
$\cC=A\sstac R H$.
Since ${\mathcal H}$ is a pure Hopf algebroid by assumption, also this coring
extension is pure.
Note that this coring extension does not correspond to any
entwining structure.

Under the additional assumption that the right ${\mathcal
H}_R$-comodule algebra $A$ is also an $L$-ring with left
$L$-linear ${\mathcal H}_R$-coaction, 
$A$ becomes an $L$-$\cC$ bicomodule, hence
one can associate to it a Morita context ${\widetilde{\mathbb M}}(A)$ as in
\equref{eq:contextmod}.
In this situation one can associate to $A$ also another Morita context
like in \cite[Remark 3.2 (1)]{BohmBrz:cleft}.
Latter one is formulated in terms of the two convolution products,
$(f,g)\mapsto \mu_A\circ (f\otimes_L g)\circ \gamma_L$, for
$f\in{\rm Hom}_L(H,A)$ and $g\in{}_L{\rm Hom}(H,A)$ on one hand,
and $(f',g')\mapsto \mu_A\circ (f'\otimes_R g')\circ \gamma_R$,
for $f'\in{\rm Hom}_R(H,A)$ and $g'\in{}_R{\rm Hom}(H,A)$ on the
other hand. These two convolution products define the convolution
algebras ${}_L{\rm Hom}_L(H,A)$ and ${}_R{\rm Hom}_R(H,A)$,
respectively, and also their bimodules ${}_L{\rm Hom}_R(H,A)$ and
${}_R{\rm Hom}_L(H,A)$. The connecting homomorphisms of the Morita
context are defined as projections of the appropriate convolution
product. 

The precise relation of the two Morita contexts in the previous paragraph
is formulated in the following lemma. 
\begin{lemma} \label{lem:BB_context}
Let ${\mathcal H}=({\mathcal H}_L,{\mathcal H}_R,S)$ be a 
pure
Hopf algebroid and $A$ be a right ${\mathcal H}_R$-comodule algebra.
Denote the associated $A$-coring $A\otimes_R H$ by $\cC$. Assume
that $A$ is also an $L$-ring and the ${\mathcal H}_R$-coaction on $A$ is
left $L$-linear. Then the Morita context $\widetilde{\mathbb M}(A)$,
associated to the $L$-$\cC$ bicomodule $A$ as in \equref{eq:contextmod}, is
isomorphic to a sub-Morita context of
\begin{equation}\label{eq:Morconv}
({}_L{\rm Hom}_L(H,A)\ ,\ {}_R{\rm Hom}_R(H,A)\ ,\ {}_L{\rm
Hom}_R(H,A)\ ,\ {}_R{\rm Hom}_L(H,A)\ ,\ \blackddiamond\ ,\
\ddiamond),
\end{equation}
where the algebra and bimodule structures are given by the respective 
convolution product and the connecting homomorphisms $\blackddiamond$ and
$\ddiamond$ are defined as projections of the convolution products.
\end{lemma}
\begin{proof}
The endomorphism algebra $T={\rm End}^\cC(A)$ can be identified with 
the subalgebra of $\cC$ (equivalently, ${\mathcal H}_R$) -coinvariants in $A$,
via $T\ni t\mapsto t(1_A)$  cf. \cite[28.4]{BrzWis:cor}. This injection
$T\hookrightarrow A$ of $L$-rings defines an injection of convolution 
algebras  
$$
\iota_1:{}_L{\rm Hom}_L(H,T)\hookrightarrow {}_L{\rm Hom}_L(H,A).
$$
The algebra of left $\cC$-colinear right ${\mathcal H}_L$-colinear
endomorphisms of $A\otimes_R H$ can be injected into the other
convolution algebra ${}_R{\rm Hom}_R(H,A)$ via the map
$$
\iota_2: {}^\cC{\rm End}^{{\mathcal H}_L}(\cC)\hookrightarrow {}_R{\rm
Hom}_R(H,A),\qquad u\mapsto \big[\ h\mapsto \big((A\stac{R}
\pi_R)\circ u\big)(1_A\stac{R}h)\ \big].
$$
For a pure Hopf algebroid ${\mathcal H}$, any right ${\mathcal H}_L$-colinear
map is also right ${\mathcal H}_R$-colinear, thus in particular right
$R$-linear.  
So we have an obvious injection
$$
\iota_3:{}_L{\rm Hom}^{{\mathcal H}_L}(H,A)\hookrightarrow {}_L{\rm Hom}_R
(H,A).
$$
Finally, by standard hom-tensor relations, we have an inclusion
$$
\iota_4:\widetilde{Q}\hookrightarrow {}_A{\rm Hom}_L(A\stac{R}H,A)\simeq
     {}_R{\rm Hom}_L(H,A).
$$
It is left to the reader as an easy exercise to check that the
four injections constructed define a morphism of Morita contexts.
\end{proof}

Following \cite{BohmBrz:cleft},
a right ${\mathcal H}_R$-comodule algebra $A$ (with $R$-ring structure
$\eta_R:R\to A$) for a pure Hopf algebroid ${\mathcal H}=({\mathcal
  H}_L,{\mathcal H}_R,$ $S)$ is called an {\em ${\mathcal H}$-cleft}
extension of its ${\mathcal H}_R$-coinvariant subalgebra $T$ if
the following conditions are satisfied.
\begin{blist}
\item $A$ is an $L$-ring (with unit morphism $\eta_L:L\to A$) and
$T$ is an $L$-subring of $A$; 
\item There exists a left $L$-linear right ${\mathcal H}_R$-colinear
  map $\lambda:H\to A$ which is invertible in the Morita context
  \eqref{eq:Morconv}, i.e. for which there exists a left $R$-linear
  right $L$-linear map ${\bar \lambda}:H\to A$ such that
$$
\lambda\ddiamond {\bar\lambda}\equiv
\mu_A\circ (\lambda\stac{R} {\bar \lambda})\circ \gamma_R=
\eta_L\circ\pi_L\quad \textrm{and}  \quad 
{\bar\lambda}\blackddiamond\lambda\equiv\mu_A\circ ({\bar
\lambda}\stac{L} {\lambda})\circ \gamma_L=\eta_R\circ\pi_R.
$$
\end{blist}
(For more details on cleft extensions by Hopf algebroids, in particular for
their characterization as crossed products, we refer
to \cite{BohmBrz:cleft}.) This definition
can be reformulated using 
the terminology of the present paper as follows. The proof is a
simple generalization of that of Proposition
\ref{prop:cleftentw}, using \cite[Lemmata 3.6 and 3.7]{BohmBrz:cleft}.
\begin{proposition}\label{prop:hgdcleft}
Let ${\mathcal H}=({\mathcal H}_L,{\mathcal H}_R,S)$ be a pure
Hopf algebroid and $A$ be a right ${\mathcal H}_R$-comodule algebra. Let
$\cC$ be the associated $A$-coring $A\sstac{R} H$, with coproduct
inherited from ${\mathcal H}_R$. $A$ is an ${\mathcal H}$-cleft
extension of its ${\mathcal H}_R$-coinvariant subalgebra if and
only if the following hold.
\begin{blist}
\item $A$ is an $L$-ring; \item $A$ (with the left $L$-module
structure in (a), the right regular
  $A$-module structure and the given ${\mathcal H}_R$-comodule
  structure) is a cleft $L$-$\cC$ bicomodule for
  the coring extension ${\mathcal H}_L$ of $\cC$.
\end{blist}
\end{proposition}
It should be emphasized that cleft extensions of algebras by a pure
Hopf algebroid are {\em not} examples of the kind discussed in
Section \ref{sec:Lcleft}. As explained, one can associate an
$R$-entwining structure -- with $R$-ring $A$ and $R$-coring
underlying ${\mathcal H}_R$ -- to a right ${\mathcal
H}_R$-comodule algebra $A$, for a Hopf algebroid 
${\mathcal H}=({\mathcal H}_L,{\mathcal H}_R,S)$.
It isn't true, however, that this $R$-entwining structure was cleft
for an ${\mathcal H}$-cleft extension. The right ${\mathcal
H}_R$-coaction in $A$ is determined by the grouplike element $1_H$
(cf. last paragraph in Section \ref{sec:Lcleft}), which is {\em
not} central in the $R$-$R$ bimodule ${\mathcal H}_R$. Hence the
right ${\mathcal H}_R$-coaction in $A$ is {\em not} left
$R$-linear in the sense of a cleft $R$-entwining structure (i.e. of  
Proposition \ref{prop:cleftLentw} (b)).

\begin{remark}\label{rem:non.pure}
Let us note that a Morita context, generalizing that in Lemma
\ref{lem:BB_context}, can be constructed also for a not necessarily pure Hopf
algebroid ${\mathcal H}$ and a right ${\mathcal H}$-comodule algebra $A$. 
Although, (since in this case the coring underlying ${\mathcal H}_L$ is not
necessarily a pure extension of the $A$-coring $\cC:={A\ot_R {\mathcal H}_R}$,
thus the functor $U:{\mathcal M}^\cC\to {\mathcal M}^{{\mathcal H}_L}$ is not
known to exist) this Morita context is not known to be of the form
\equref{eq:contextmod}. 

Consider a (not necessarily pure) Hopf algebroid ${\mathcal H}$, with
constituent right $R$-bialgebroid ${\mathcal H}_R$ and left $L$-bialgebroid
${\mathcal H}_L$. 
As it is described in \cite[{\em Corrigendum} Definition 1]{BohmBrz:cleft},
the right notion of an ${\mathcal H}$-comodule is a right ${\mathcal
  H}_R$-comodule and right ${\mathcal H}_L$-comodule $M$, such that 
the ${\mathcal H}_R$-coaction is an ${\mathcal H}_L$-comodule map and the
${\mathcal H}_L$-coaction is an ${\mathcal H}_R$-comodule map. Comodules of
${\mathcal H}$ defined in this way constitute a monoidal category
$\Mm^{\mathcal H}$, with strict monoidal forgetful functors $\Mm^{\mathcal
  H}\to \Mm^{{\mathcal H}_R}$ and $\Mm^{\mathcal H}\to \Mm^{{\mathcal H}_L}$,
see \cite[{\em Corrigendum} Theorem 6]{BohmBrz:cleft}. 
(If ${\mathcal H}$ is a pure Hopf algebroid, then both of these forgetful
functors become isomorphisms, cf. \cite[{\em Corrigendum} Theorem
  4(3)]{BohmBrz:cleft}.) 
One defines an ${\mathcal H}$-comodule algebra $A$ as a monoid (algebra
object) in $\Mm^{\mathcal H}$. Considering the category $\Mm^{\mathcal H}_A$
of right $A$-modules in $\Mm^{\mathcal H}$, there is a forgetful functor
$U':\Mm^{\mathcal H}_A \to \Mm^{\mathcal H}$. (If ${\mathcal H}$ is a pure Hopf
algebroid, then $U'$ differs from the functor $U:= - \Box_{\cC}\, \cC:\Mm^\cC\to
\Mm^{{\mathcal H}_L}$ by trivial isomorphisms. However, in the general case
$\Mm^{\mathcal H}_A$ is not known to be isomorphic to the category of
comodules for any coring, and $U'$ is not known to be given by a cotensor
product.)  

If $A$ is in addition an $L$-ring, with left $L$-linear ${\mathcal
  H}_R$-coaction, then 
there is another functor $V':=\mathrm{Hom}^{\mathcal H}_A(A,-)\ot_L
H:\Mm^{\mathcal H}_A \to \Mm^{\mathcal H}$. (If ${\mathcal H}$ is a pure Hopf
algebroid, then $V'$ differs from the functor $V:=\mathrm{Hom}^\cC(A,-)\ot_L
H:\Mm^\cC \to \Mm^{{\mathcal H}_L}$ by trivial isomorphisms.)
Both functors $U'$ and $V'$ determine a Morita context as in
\eqref{eq:catMor}. Similar computations used to prove \prref{contextmod} and
Lemma \ref{lem:BB_context} verify that this Morita context is isomorphic to 
\begin{equation}\eqlabel{eq:BB}
\big(\ {}_L \Hom_L(H,T)\ ,\ S\ ,\ {}_L \Hom^{\mathcal H}(H,A)\ ,\ {}_{L^{op}}
\Hom^{\mathcal H}(H^{tw},A)\ ,\ \blackddiamond\ ,\ \ddiamond\ \big).
\end{equation}
Here $T$ is the ${\mathcal H}_R$-coinvariant subalgebra of $A$, $S$ is a
subalgebra of the convolution algebra ${}_R\Hom_R(H,A)$ and $H^{tw}$ is a right
${\mathcal H}$-comodule with underlying $k$-module $H$ and coactions obtained by
twisting the coproducts (regarded as left coactions) with the antipode,
cf. \cite[{\em Corrigendum} Remark 2]{BohmBrz:cleft}. The connecting maps
$\blackddiamond$ and $\ddiamond$ are restrictions of the convolution product in
Lemma \ref{lem:BB_context}.
For more details see the arXiv version of \cite{Bohm:rev}.

Definition of an ${\mathcal H}$-cleft extension $T\subseteq A$
in \cite[\emph{Corrigendum}]{BohmBrz:cleft} is equivalent to the existence of
invertible elements in the Morita context \equref{eq:BB}, hence in the
isomorphic Morita context determined by the functors $U'$ and $V'$. 
\end{remark}

In light of Remark \ref{rem:non.pure}, in particular the Strong Structure
Theorem in Corollary \ref{cor:cleft.str.thm} (2) has no message about cleft
extensions by a Hopf algebroid ${\mathcal H}$ (unless ${\mathcal H}$ is
assumed to be pure). Instead, making use of the Morita context \equref{eq:BB},
one proves the following. 

\begin{theorem}\thlabel{thm:nonpure.str.str.thm}
Let ${\mathcal H}$ be any Hopf algebroid and $T\subseteq A$ be an ${\mathcal
  H}$-cleft extension (in the sense of
\cite[\emph{Corrigendum}]{BohmBrz:cleft}). Then there is an equivalence
functor $-\ot_T A:\cM_T \to \cM^{\mathcal H}_A$. 
\end{theorem}

\begin{proof}
Consider mutually inverse elements $j\in {}_L\mathrm{Hom}^{\mathcal H}(H,A)$ and
${\widetilde j}\in {}_{L^{op}} \mathrm{Hom}^{\mathcal H}(H^{tw},A)$ in the Morita 
context \equref{eq:BB}.
For the coaction of the constituent right $R$-bialgebroid ${\mathcal H}_R$,
and the coaction of the constituent left $L$-bialgebroid ${\mathcal H}_L$ on
$M\in \cM^{\mathcal H}_A$, use the Sweedler's index notations $m\mapsto
m^{[0]}\ot_R m^{[1]}$ 
(with upper indeces) and $m\mapsto m_{[0]}\ot_L m_{[1]}$ (with lower indeces),
respectively, where in both cases implicit summation is understood. In
particular, denote the coproducts in ${\mathcal H}_R$ and ${\mathcal H}_L$ by
$h\mapsto h^{(1)}\ot_R h^{(2)}$ and $h\mapsto h_{(1)}\ot_L h_{(2)}$, respectively.

For any right $T$-module $W$, $W\ot_T A$ is an object in $\cM^{\mathcal H}_A$
via the $A$-action and the ${\mathcal H}$-coactions on the second factor. The
resulting functor $-\ot_T A:\cM_T \to \cM^{\mathcal H}_A$ is left adjoint of
the ${\mathcal H}_R$-coinvariants functor $(-)^{co{\mathcal H}_R}:
\cM^{\mathcal H}_A  \to \cM_T$, where the $T$-action on $M^{co{\mathcal H}_R}$
is induced by the $A$-action on $M\in \cM^{\mathcal H}_A$. The counit and the
unit of the adjunction are given by 
$$
\begin{array}{lll}
c_M:M^{co{\mathcal H}_R}\stac  T A\to M,\qquad &n \stac  T a \mapsto na,\qquad
  &\textrm{for } M \in \cM^{\mathcal H}_A;\\
e_W:W \to (W\stac  T A)^{co{\mathcal H}_R},\qquad &w \mapsto w\stac  T 1_A,\qquad
  &\textrm{for } W\in \cM_T.
\end{array}
$$
The map $c_M$ is a morphism in $\cM^{\mathcal H}_A$ by \cite[{\em Corrigendum}
  Proposition 3]{BohmBrz:cleft}.
We prove that $-\ot_T A:\cM_T \to \cM^{\mathcal H}_A$ is an equivalence by
constructing the inverse of the above natural transformations.

The inverse of $c_M$ is given by the map
$$
c_M^{-1}(m)
= m^{[0]} {\widetilde j}({m^{[1]}}_{(1)}) \stac  T j(m{^{[1]}}_{(2)})
= {m_{[0]}}^{[0]} {\widetilde j}({m_{[0]}}^{[1]}) \stac  T j(m_{[1]}),
$$
where the second equality follows by the right ${\mathcal H}_L$-colinearity of
the ${\mathcal H}_R$-coaction on $M$, cf. \cite[{\em Corrigendum} Definition
  1]{BohmBrz:cleft}. 
The map $c_M^{-1}:M \to M \ot_T A$ is well defined by the $R$-, and $L$-module
map properties of ${\widetilde j}$ and $j$, and since $T$ is an $L$-ring.
For any $m\in M$, the element $m^{[0]} {\widetilde j}({m^{[1]}})$ is checked
to be ${\mathcal H}_R$-coinvariant, using 
the right ${\mathcal H}_R$-colinearity of the $A$-action on $M$ and the right
${\mathcal H}_R$-colinearity of
${\widetilde j}$, the right ${\mathcal H}_L$-colinearity of the coproduct in
${\mathcal H}_R$, one of the 
antipode axioms in a Hopf algebroid, an anti-isomorphism between the base
algebras $L$ and $R$ of both constituent bialgebroids, together with 
the right $R$-linearity of ${\widetilde j}$ and the counitality of the
coproduct in ${\mathcal H}_L$.
Hence the range of $c_M^{-1}$ is in $M^{co{\mathcal H}_R}\ot_T A$, as
needed. 
For $m\in M$, and $n\ot_T a\in M^{co{\mathcal H}_R}\ot_T A$,  
\begin{eqnarray*}
(c_M\circ c_M^{-1})(m)
&=& m^{[0]} {\widetilde j}({m^{[1]}}_{(1)}) j({m^{[1]}}_{(2)}) 
= m;\\
( c_M^{-1} \circ c_M)(n\stac  T a)
&=& n {a_{[0]}}^{[0]} {\widetilde j}({a_{[0]}}^{[1]}) \stac T j(a_{[1]})\\
&=&n\stac T {a_{[0]}}^{[0]} {\widetilde j}({a_{[0]}}^{[1]})j(a_{[1]})
=n\stac T {a}^{[0]} {\widetilde j}({a^{[1]}}_{(1)})j({a^{[1]}}_{(2)})
= n\stac  T a.
\end{eqnarray*}
In the first equality in the last line we used that 
$a^{[0]}{\widetilde j}(a^{[1]})$ belongs to $T$, for any $a\in A$.
The second equality in the last line follows by the right ${\mathcal
  H}_L$-colinearity of the ${\mathcal H}_R$-coaction on $M$.
In the last equality of both computations we used that ${\widetilde j}$ and
$j$ are mutual inverses in the Morita context \equref{eq:BB}, together with
the counitality of the ${\mathcal H}_R$-coaction.  

The inverse of $e_W$ is given by 
$$
e_W^{-1}( \sum_i w_i\stac  T a_i) 
= \sum_i w_i \, a_i^{[0]}{\widetilde j}(a_i^{[1]})j(1_H). 
$$
Since
$a^{[0]}{\widetilde j}(a^{[1]})$ belongs to $T$, for any $a\in A$, and 
by the ${\mathcal H}_R$-colinearity of $j$ and the unitality of the coproduct
in ${\mathcal H}_R$, the element $j(1_H)$ also belongs to $T$, the expression
of $e_W^{-1}$ is meaningful. 
Obviously, $(e_W^{-1}\circ e_W)(w)=w$, for all $w\in W$. To check that
$e_W^{-1}$ is also the right inverse of $e_W$, note that for $\sum_i w_i\ot_T
a_i \in (W\ot_T A)^{co{\mathcal H}_R}$ the identity 
$\sum_i w_i \ot_ T a_i^{[0]}\ot_ R a_i^{[1]}
= \sum_i w_i \ot_ T a_i \ot_ R 1_H$ holds, 
which implies 
\begin{eqnarray*}
(e_W\circ e_W^{-1})(\sum_i w_i\stac  T a_i)
&=& \sum_i w_i  a_i^{[0]}{\widetilde j}(a_i^{[1]}) j(1_H) \stac  T 1_H \\
&=& \sum_i w_i \stac  T a_i^{[0]}{\widetilde j}(a_i^{[1]}) j(1_H)
= \sum_i w_i \stac  T a_i.
\end{eqnarray*} 
\end{proof}

The total algebra $H$ of a Hopf algebroid ${\mathcal H}$ is an ${\mathcal
  H}$-cleft extension of the base algebra $L$ of the constituent left
bialgebroid ${\mathcal H}_L$,  
via the source map of ${\mathcal H}_L$. Indeed,
mutually inverse elements in the corresponding Morita context \equref{eq:BB}
are provided by the identity map in  ${}_L\mathrm{Hom}^{\mathcal H}(H,H)$ and
the antipode in ${}_{L^{op}} \mathrm{Hom}^{\mathcal H}(H^{tw},H)$. 
Thus 
\thref{thm:nonpure.str.str.thm} implies, in particular, that for any Hopf
algebroid ${\mathcal H}$, there is an equivalence functor $-\ot_L H: \cM_L \to
\cM^{\mathcal H}_H$. This yields a {\em corrected version} of the Fundamental
Theorem of Hopf modules, \cite[Theorem 4.2]{Bohm:hgdint}.
Regrettably, the proof of the journal version of \cite[Theorem
  4.2]{Bohm:hgdint} turned out to be incorrect (when checking that the
to-be-inverse of the counit of the adjunction has the required range, 
some ill-defined maps are used). 
Hence in the published (stronger) form, \cite[Theorem 4.2]{Bohm:hgdint} is not 
justified without some further assumption, see the {\em Corrigendum} and the
arXiv version.

\section*{Acknowledgements} 
GB is grateful to the members of the Department of Mathematics,
Faculteit Ingenieurswetenschappen VUB, for a very pleasant visit in September
2005. She thanks Korn\'el Szlach\'anyi for valuable discussions about
tripleability theorems. Her work is supported by the Hungarian Scientific
Research Fund OTKA T 043 159 and the Bolyai J\'anos Fellowship.

JV would like to thank Gabi B\"ohm for warm hospitality during
his stay in Budapest where this work initiated.

\end{document}